\documentclass{article}
\usepackage{graphicx} 
\usepackage{amsmath,amssymb,amsfonts,amscd,amsxtra,amsthm}
\usepackage{dsfont}
\usepackage{latexsym}
\usepackage{xcolor}
\usepackage[title]{appendix}
\usepackage{floatrow}
\usepackage{graphicx}
\usepackage{subcaption} 
\usepackage{booktabs}
\usepackage{bbm}
\usepackage{float}
\usepackage{overpic}

\usepackage{tikz-cd}
\usepackage{tikz}
\usetikzlibrary{arrows.meta,calc}
\usepackage{algorithmic}
\usepackage[ruled,linesnumbered,noend]{algorithm2e}
\usepackage{zymacros}
\usepackage{appendix}
\usepackage{epstopdf}
\usepackage[top=1in, bottom=1in, left=1in, right=1in]{geometry}
\renewcommand{\theequation}{\thesection.\arabic{equation}}

\renewcommand\appendix{\par
  \setcounter{section}{0}
  \setcounter{subsection}{0}
  \setcounter{figure}{0}
  \setcounter{table}{0}
  \renewcommand\thesection{Appendix \Alph{section}}
  \renewcommand\theequation{\Alph{section}.\arabic{equation}}
  \renewcommand\thefigure{\Alph{section}.\arabic{figure}}
  \renewcommand\thetable{\Alph{section}.\arabic{table}}
  \renewcommand\thethm{\Alph{section}.\arabic{thm}}
}

\numberwithin{equation}{section}

\newcommand{\bra}[1]{\left(#1\right)}
\usepackage{authblk}

\title{Elastic Wave Imaging via Kirchhoff Migration: Resolution Analysis and Inclusion Identification}
\author{
Zetao Fei\thanks{\footnotesize CMAP, Ecole Polytechnique, Institut Polytechnique de Paris, 91120 Palaiseau, France (zetao.fei@polytechnique.edu).}\;
	and Josselin Garnier\thanks{\footnotesize 
	CMAP, Ecole Polytechnique, Institut Polytechnique de Paris, 91120 Palaiseau, France (josselin.garnier@polytechnique.edu). This work was partially supported by Agence Nationale de la Recherche (Grant No. ANR-23-CE30-0021).}}

\begin{document}

\maketitle
\begin{abstract}
We consider the inverse problem of determining the location,
orientation, shape, and material parameters of a small elastic
inclusion embedded in a homogeneous medium from full-aperture
measurements of the scattered field. Using a small-volume expansion,
we derive explicit leading-order expressions for four Kirchhoff
migration imaging functionals, denoted by PP, PS, SP, and SS,
corresponding to the four combinations of pressure (P) and shear (S)
waves used for receiver- and source-side backpropagation. The
resulting point-spread functions quantitatively characterize the
resolution of each mode and show that the SS mode achieves the highest
spatial resolution. For an elliptic inclusion, we show that the PP
mode uniquely determines the elastic moment tensor (EMT) within the
leading-order asymptotic model, thereby enabling the reconstruction of
the material and geometric parameters of the inclusion. For a general
small scatterer, the reconstruction may instead be interpreted as an
effective elliptic representation determined by its leading-order EMT.
For a disk, the parameter-identification problem is underdetermined;
nevertheless, explicit upper and lower bounds on its volume are
derived. We further analyze the model error arising from neglected
cross-mode contributions when the distance from the center of the
measurement configuration is not sufficiently large relative to the
elastic wavelengths. Numerical experiments support the theoretical
predictions and demonstrate the effectiveness of the proposed imaging
and reconstruction methods.
\end{abstract}

\section{Introduction}

Inverse scattering problems for elastic waves have been extensively studied due to their broad applications in nondestructive testing \cite{INESSILVA2023101155}, medical imaging \cite{greenleaf2003selected}, and geophysical exploration \cite{aki2002quantitative}. Significant progress has been made in the mathematical theory of elastic inverse scattering, including uniqueness and stability; see, for instance, \cite{hahner1993uniqueness,nakamura1994global,elschner2010uniqueness,alessandrini2014stable,barcelo2018uniqueness,zhai2026determination} and the references therein.

When the characteristic size of the inclusion is small compared with both the compressional and shear wavelengths, the scattered field admits a small-volume asymptotic expansion \cite{3d546238-6f8a-3f6d-b0b4-3883234b540d}. Based on this expansion, the inverse problem of reconstructing a small inclusion for the time-harmonic elasticity equations can be formulated as a nonlinear least-squares optimization problem. However, convergence of the resulting iterative reconstruction to the true solution is generally not guaranteed \cite{beretta2010algorithm}. Over the past decades, a variety of imaging methods for small elastic inclusions have been developed and analyzed, including MUSIC-type methods \cite{doi:10.1137/070696076}, time-reversal methods \cite{ammari2011transient}, and topological derivative-based methods \cite{guzina2004topological,ammari2013localization}. As a simple, computationally efficient, and robust imaging technique, Kirchhoff migration has also attracted considerable attention in wave-based inverse problems \cite{claerbout1985imaging,bleistein2001mathematics}. Owing to the coexistence of compressional (P) and shear (S) waves, elastic wave imaging involves multiple wave modes with different propagation speeds, making both the analysis and the imaging functional substantially more complicated than in the acoustic case. Consequently, despite its practical importance, Kirchhoff migration for elastic waves has received comparatively less
mathematical analysis than its acoustic counterpart. We refer to \cite{kuo1983kirchhoff,du2008elastic,sun2011prestack,du2012multicomponent,rabinovich2024kirchhoff} and the references therein for the development of Kirchhoff migration methods for elastic waves. In particular, a quantitative analysis of Kirchhoff migration imaging functionals associated with different polarization modes, together with their implications for the quantitative characterization of elastic inclusions, remains largely unexplored.

Besides inclusion localization and resolution analysis, another fundamental inverse problem in the study of composite materials is the estimation of the volume (or volume fraction) of embedded inclusions from boundary measurements. In \cite{doi:10.1137/S0036141001388944,ALESSANDRINI2007288,Alessandrini_2005}, upper and lower bounds for the inclusion volume were established in terms of the mechanical work associated with prescribed boundary tractions, and the resulting estimates were validated numerically. Under the assumption that the bulk and shear moduli of both constituent materials are known, \cite{MILTON2012193} derived bounds on the volume fraction. For inclusions with constant shear modulus, \cite{Thaler_2014} showed that the inclusion volume can be determined exactly from boundary measurements. In the asymptotic regime of small inclusions, the methods based on small-volume expansions presented in \cite{kang2003identification,ammari2004reconstruction,ammari2007polarization} represent an arbitrary inclusion by an equivalent ellipse. When the Lam\'e parameters are known, the geometric parameters of the equivalent ellipse can be reconstructed. If the Lam\'e parameters are unknown, however, the reconstruction is carried out only under the assumption that the inclusion is either a cavity (vanishing Lam\'e parameters) or a rigid inclusion (infinite shear modulus). Further developments on estimating the geometric parameters through the elastic moment tensor (EMT) can be found, for example, in \cite{kang2007numerical,capdeboscq2008improved}.

In this paper, we investigate four Kirchhoff migration imaging
functionals, referred to as the PP, PS, SP, and SS modes, corresponding
to the four possible combinations of pressure (P) and shear (S) waves
used for receiver-side and source-side backpropagation; see
\eqref{eq:imaging_func} for their precise definitions. Our analysis
reveals a fundamental distinction between their roles: the SS mode
provides the highest spatial resolution for inclusion localization,
whereas the PP mode is the only individual mode containing sufficient
independent information for quantitative identification through the
elastic moment tensor (EMT). Specifically, we
(i) derive explicit leading-order asymptotic expressions and
point-spread functions for all four modes;
(ii) compare their resolution properties and show that the SS mode
achieves the highest spatial resolution;
(iii) prove that, for an elliptic inclusion, the PP mode uniquely
determines the EMT within the leading-order asymptotic model;
(iv) develop an efficient procedure for reconstructing the material
and geometric parameters of the ellipse;
(v) show that the corresponding identification problem for a
disk-shaped inclusion is underdetermined and derive explicit upper and
lower bounds for its volume; and
(vi) analyze the model error associated with cross-mode contributions
that may become significant when
condition~\eqref{eq:large-frequency-distance} is not sufficiently
satisfied. Numerical experiments support the theoretical predictions
and demonstrate the effectiveness of the proposed imaging and
reconstruction procedures.

The remainder of the paper is organized as follows. In Section~\ref{sec:setup}, we introduce the mathematical formulation of the problem and briefly review the EMT. Section~\ref{sec:psf_resolution} is devoted to the asymptotic analysis of the imaging functionals and their resolution properties for the four polarization modes. In Section~\ref{sec:identification}, we study the identification problem for elliptic and disk-shaped inclusions. Numerical experiments are presented in Section~\ref{sec:numerics}, and concluding remarks are given in Section~\ref{sec:conclusion}.

\section{Problem formulation}\label{sec:setup}
In this section, we present the mathematical preliminaries for the model considered in this paper. Section~\ref{subsec:background} introduces the problem setup and the imaging configuration. Section~\ref{subsec:intro_EMT} provides a brief review of the EMT, together with the essential results used throughout the paper.
\subsection{Basic setup for elastic wave equations with a small inclusion}\label{subsec:background}
We consider an elastic medium in $\sR^2$. Let the background Lam\'e parameters be $(\lambda,\mu)$ satisfying 
\begin{align}\label{eq:lame_constraint}
    \mu>0, \quad \lambda+\mu>0,
\end{align}
and the background volume density be $\rho>0$. The physical significance of the constraints on the Lam\'e parameters (\ref{eq:lame_constraint}) is discussed in \cite{landau2012theory}.
Let $\mGamma^{\omega}(\vx,\vy):=\mGamma^{\omega}(\vx-\vy)$ be the fundamental solution to the homogeneous time-harmonic elastic wave equation in $\sR^2$, i.e., the solution to 

\[
\renewcommand{\arraystretch}{1.5}
\begin{aligned}
    \left\{
    \begin{array}{l}
     (\fL^{\lambda,\mu}+\omega^2\rho)\mGamma^{\omega}(\vx,\vy)= \delta(\vx-\vy) \mI_2     \quad\text{in } \sR^2,  \\
    \mGamma^{\omega} \text{ satisfies the Sommerfeld--Kupradze radiation condition},\\
    \end{array}
    \right.
\end{aligned}
\]
where $\mI_2$ is the $2\times2$ identity matrix. Here the operator $\fL^{\lambda,\mu}$ is defined by 
\[
    \fL^{\lambda,\mu}[\vu] := \mu\Delta\vu+(\lambda+\mu)\nabla\nabla\cdot\vu.
\]
We denote the background wavenumbers $\kappa_P$ and $\kappa_S$ by
\[
    \kappa_P = \frac{\omega}{c_P}, \quad \kappa_S = \frac{\omega}{c_S},
\]
where $c_P$ and $c_S$ are the compressional and shear wave speeds, respectively:
\[
    c_P = \sqrt{\frac{\lambda+2\mu}{\rho}}, \quad c_S=\sqrt{\frac{\mu}{\rho}}.
\]
Throughout the paper, we write $\vx=(x_1,x_2)$, use $\partial_j$ to
denote partial differentiation with respect to $x_j$ $(j=1,2)$, and
denote by $\ve_j$ the unit vector in the $x_j$-direction.

In two dimensions, upon identifying the fundamental solution with a
function of the relative variable, its components are given by
\begin{align}
    \Gamma_{ij}^\omega(\vx)
    =
    -\frac{i}{4\mu}\delta_{ij}H_{0}^{(1)}(\kappa_S|\vx|)
    +\frac{i}{4\omega^2\rho}
    \partial_i\partial_j
    \bigl(
        H_{0}^{(1)}(\kappa_P|\vx|)
        -
        H_{0}^{(1)}(\kappa_S|\vx|)
    \bigr),
\end{align}
where $\delta_{ij}$ is the Kronecker delta and $H^{(1)}_0$ denotes the
Hankel function of the first kind and order zero.

We consider an active imaging configuration with prescribed source
profiles and full-aperture measurements. A collocated source--receiver
array $\{\vx_k\}_{k=1}^N$ is distributed along the circular measurement
curve $\partial\Omega$ of radius $R$, centered at the origin $O$. We choose an arbitrary reference direction through the center $O$ of
the circular measurement configuration as the positive $x_1$-axis,
and define the positive $x_2$-axis by a counterclockwise rotation of
$\pi/2$. The resulting Cartesian coordinate system is fixed throughout
the analysis and is referred to as the laboratory coordinate system. We
take $R$ as the reference length scale, with $R\sim\mathcal{O}(1)$.

Relative to this scale, we consider a small elastic inclusion $D$ of
the form
\begin{equation}\label{eq:small_inclusion}
    D=\vz+\varepsilon B,
\end{equation}
where $B\subset\mathbb{R}^2$ is an
ellipse centered at the origin whose major and minor axes are aligned
with the laboratory coordinate axes, $\vz$ denotes the center of the inclusion in the laboratory coordinate
system, and $\varepsilon\ll1$ is a dimensionless scaling parameter.
Thus, if
\[
    L_B:=\operatorname{diam}(B),
\]
then the characteristic diameter of $D$ is $\varepsilon L_B$.

When the inclusion is an ellipse having orientation angle $\phi\in[0,2\pi)$, we explicitly distinguish between the
laboratory and principal-axis coordinate systems. Let
\begin{align}\label{eq:def_R_phi}
    \mR^\phi
    :=
    \begin{pmatrix}
        \cos\phi & -\sin\phi\\
        \sin\phi &  \cos\phi
    \end{pmatrix}
\end{align}
be the rotation matrix. For the orientation angle $\phi$, we define the rotated ellipse by
\[
    B_\phi:=\mR^\phi(B),
\]
and the corresponding physical inclusion is represented as
\begin{equation}\label{eq:rotated_inclusion}
    D_\phi=\vz+\varepsilon B_\phi
    =\vz+\varepsilon\mR^\phi(B).
\end{equation}
For the elliptic inclusion $D_\phi$, we introduce a principal-axis coordinate
system centered at $\vz$, whose first and second coordinate axes are
aligned with the major and minor axes of the ellipse, respectively.
If a point has coordinates $\vx$ in the laboratory coordinate system,
then its coordinate vector $\boldsymbol{\xi}=(\xi_1,\xi_2)^\top$ in
the principal-axis coordinate system is given by
\begin{equation}\label{eq:principal_coordinates}
    \boldsymbol{\xi}
    =
    (\mR^\phi)^\top(\vx-\vz),
\end{equation}
or, equivalently,
\[
    \vx=\vz+\mR^\phi\boldsymbol{\xi}.
\]
For a
disk, the choice of $\phi$ is immaterial. We emphasize that, although the theoretical results may be expressed in the principal-axis coordinate system for analytical convenience, the imaging configuration, the definitions of the imaging functionals, and the sampling procedure are all formulated in the laboratory coordinate system, which represents the physical coordinate system used in practice.

To describe the frequency scaling, we set
\begin{equation}\label{eq:frequency_scaling}
    \omega
    =
    \varepsilon^{-1/2}\widetilde{\omega},
\end{equation}
where $\widetilde{\omega}>0$ is independent of $\varepsilon$
and has the same physical dimension as $\omega$.
We keep the inclusion center $\vz\neq\mathbf{0}$ and the
measurement radius $R$ independent of $\varepsilon$. In addition, we impose the geometric separation
$\delta:=|\vz|/R\ll1$, independently of $\varepsilon$.
The resulting scale-separation assumptions are summarized as
\begin{subequations}\label{eq:scaling}
\begin{align}
    \operatorname{diam}(D_\phi)
    &\ll
    \frac{2\pi c_\alpha}{\omega},
    && \alpha=P,S,
    \label{eq:small-volume-regime}
    \\
    \frac{2\pi c_\alpha}{\omega}
    &\ll
    |\vz|,
    && \alpha=P,S,
    \label{eq:large-frequency-distance}
    \\
    |\vz|
    &\ll R.
    \label{eq:geometric-separation}
\end{align}
\end{subequations}
Here, \eqref{eq:small-volume-regime} and
\eqref{eq:large-frequency-distance} hold asymptotically as
$\varepsilon\to0$, while \eqref{eq:geometric-separation}
is an independent geometric assumption.

Condition~\eqref{eq:small-volume-regime} places the inclusion in the
small-volume regime for both compressional and shear waves and
justifies the small-volume asymptotic expansion of the scattered
field. Condition~\eqref{eq:large-frequency-distance} is an additional asymptotic assumption used to derive explicit expressions for the imaging functionals. When this condition is not sufficiently satisfied, the actual imaging profiles may differ significantly from their asymptotic predictions, resulting in a non-negligible model error. A qualitative analysis of this error is
presented in Section~\ref{subsec:model_error}.
Condition~\eqref{eq:geometric-separation} ensures that the inclusion is
well separated from the measurement curve and, together with the
preceding wavelength assumptions, justifies the far-field asymptotic
expansion of the background Green tensor.

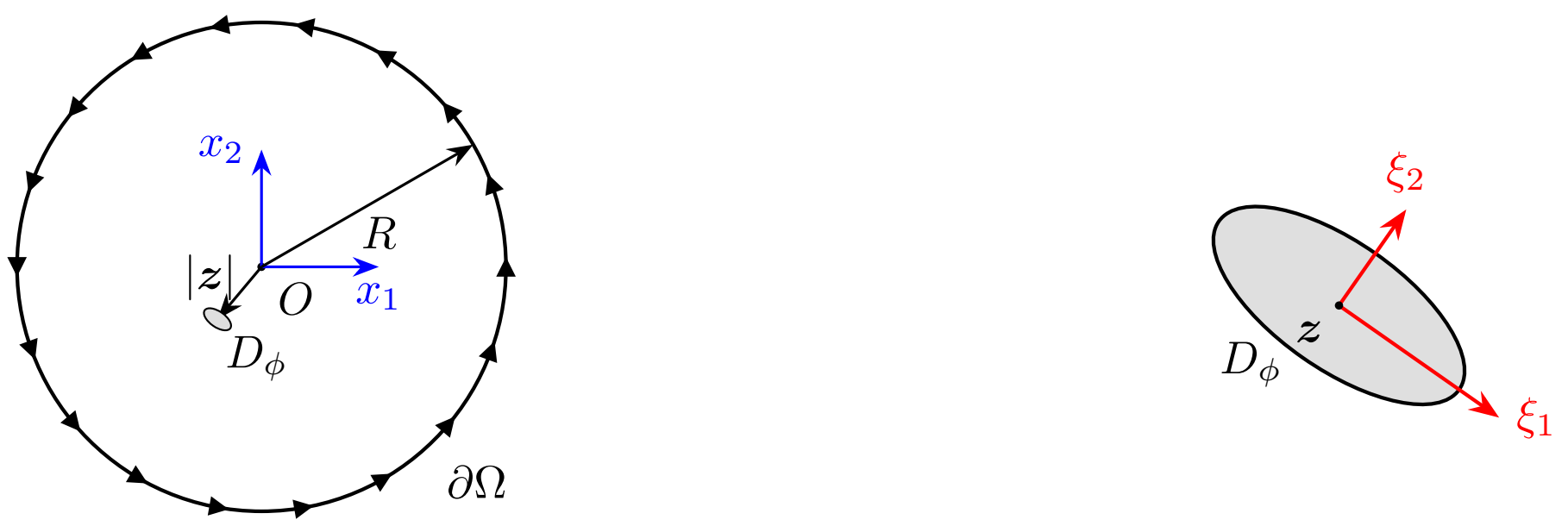
\begin{figure}[ht]
    \centering

\begin{subfigure}[t]{0.48\textwidth}
    \centering
    \begin{tikzpicture}[scale=0.75,>=Stealth]

        \def\R{2.5}
        \def\r{0.7}
        \def\alpha{-35}

        \coordinate (O) at (0,0);
        \coordinate (Z) at ({\r*cos(230)},{\r*sin(230)});

        \draw[thick] (O) circle (\R);

        \fill (O) circle (1.2pt);
        \node[below right] at (O) {$O$};

        \draw[->,blue,semithick]
            (O) -- (1.2,0)
            node[below] {$x_1$};

        \draw[->,blue,semithick]
            (O) -- (0,1.2)
            node[left] {$x_2$};

        \draw[->,semithick]
            (O) -- ({\R*cos(30)},{\R*sin(30)});
        \node at (1.2,0.35) {$R$};

        \draw[->,semithick] (O) -- (Z);
        \node[left] at ($(O)!0.2!(Z)$) {$|\vz|$};

        \foreach \angle in {0,20,...,340} {
            \begin{scope}[rotate=\angle]
                \fill
                    (\R,0) ++(0,0.1)
                    -- ++(-0.1,-0.2)
                    -- ++(0.2,0)
                    -- cycle;
            \end{scope}
        }

        \begin{scope}[shift={(Z)},rotate=\alpha]
            \filldraw[
                    fill=gray!25,
                ] (0,0) ellipse (0.16 and 0.08);
        \end{scope}

        \node at ($(Z)+(0.4,-0.4)$) {$D_{\phi}$};
        \node at (2.2,-2.2) {$\partial\Omega$};

    \end{tikzpicture}
    \caption{Imaging configuration in the laboratory coordinate
    system $(x_1,x_2)$.}
    \label{fig:imaging_geometry}
\end{subfigure}
    \hfill
    \begin{subfigure}[t]{0.48\textwidth}
    \centering
    \raisebox{0.5cm}{%
        \begin{tikzpicture}[scale=0.75,>=Stealth]

            \def\alpha{-35}
            \coordinate (Z) at (0,0);

            \begin{scope}[shift={(Z)},rotate=\alpha]

                \filldraw[
                    fill=gray!25,
                    draw=black,
                    thick
                ] (0,0) ellipse (1.5 and 0.65);

                \draw[->,red,thick]
                    (0,0) -- (2,0)
                    node[right] {$\xi_1$};

                \draw[->,red,thick]
                    (0,0) -- (0,1.2)
                    node[above] {$\xi_2$};

            \end{scope}

            \fill (Z) circle (1.2pt);
            \node[below left] at (Z) {$\vz$};

            \node at (-0.9,-0.6) {$D_\phi$};

        \end{tikzpicture}%
    }
    \caption{Principal-axis coordinate system $(\xi_1,\xi_2)$
    associated with the elliptic inclusion.}
    \label{fig:principal_axis_coordinates}
\end{subfigure}

    \caption{Geometry of the imaging problem. In panel
    (\subref{fig:imaging_geometry}): the black triangles denote collocated
    sources and receivers distributed along the circular measurement
    curve $\partial\Omega$ of radius $R$, centered at the origin $O$.
    The small inclusion $D_\phi$ is centered at $\vz$, with $|\vz|\ll R$.
    In panel (\subref{fig:principal_axis_coordinates}): the principal axes
    $(\xi_1,\xi_2)$ of the elliptic inclusion are rotated by an angle
    $\phi$ relative to the laboratory coordinate system $(x_1,x_2)$.}
    \label{fig:imaging_config}
\end{figure}

These conditions should be understood as scale-separation assumptions
defining the asymptotic regime, rather than as intrinsic physical
restrictions on the parameters of the problem. The resulting
configuration is representative of standard active imaging settings
arising in nondestructive testing and related inverse problems. A
schematic illustration of the imaging configuration is provided in
Figure~\ref{fig:imaging_config}.

The inclusion is assumed to have Lam\'e parameters $(\lambda',\mu')$ satisfying
\begin{align}\label{eq:constraint_lame_inclu}
\mu' > 0, \qquad
\lambda' + \mu' > 0, \qquad
(\lambda-\lambda')(\mu-\mu') > 0.
\end{align}
We denote the density in $D_\phi$ by $\rho'$. The scattering problem is governed by
\begin{align}\label{eq:scattering_problem}
\renewcommand{\arraystretch}{1.5}
\left\{
\begin{aligned}
&(\fL^{\lambda,\mu} + \omega^2 \rho)\vu_\varepsilon = 0
\quad \text{in } \mathbb{R}^2 \setminus \overline{D_\phi}, \\
&(\fL^{\lambda',\mu'} + \omega^2 \rho')\vu_\varepsilon = 0
\quad \text{in } D_\phi, \\
&\vu_\varepsilon|_- = \vu_\varepsilon|_+
\quad \text{on } \partial D_\phi, \\
&\frac{\partial \vu_\varepsilon}{\partial \vnu'}\Big|_-
=
\frac{\partial \vu_\varepsilon}{\partial \vnu}\Big|_+
\quad \text{on } \partial D_\phi, \\
&\vu_\varepsilon(\vx)-\vu_{\mathrm{inc}}(\vx)
\ \text{satisfies the Sommerfeld--Kupradze radiation condition}.
\end{aligned}
\right.
\end{align}
Here $\vu_{\mathrm{inc}}$ denotes the incident field, $|_-$ and $|_+$ denote the limits from the interior and exterior of the inclusion to its boundary, respectively. The conormal derivative associated with $(\lambda,\mu)$ is defined by
\[
\frac{\partial \vu_\varepsilon}{\partial \vnu}
=
\lambda (\nabla \cdot \vu_\varepsilon)\vn
+
\mu (\nabla \vu_\varepsilon + \nabla \vu_\varepsilon^\top)\vn,
\]
where $\nabla \vu_\varepsilon = (\partial_j u_{\varepsilon,i})_{i,j=1,2}$, and $\vn$ is the outward unit normal vector on $\partial D_\phi$. The conormal derivative associated with $(\lambda',\mu')$ is defined analogously.

For the $k$-th source located at $\vx_k$, we excite the medium with polarization $\ve_j$, $j=1,2$, and the incident field is given by
\[
\vu_{\mathrm{inc}}^{(j)}(\vx,\vx_k)
=
-\mGamma^\omega(\vx,\vx_k)\ve_j.
\]
The corresponding total field is measured at all receivers on $\partial\Omega$. The explicit expression for the scattered field $\vu_{\mathrm{scat}}^{(j)}$ is given in the next subsection, see (\ref{eq:scattered_expand}).

Assuming that the source--receiver array is sufficiently dense so that
the discrete sums can be approximated by boundary integrals, we define
the Kirchhoff migration imaging functionals
\begin{align}\label{eq:imaging_func}
    \fI_{\alpha,\beta}(\vz^s)
    =
    \sum_{j=1}^2
    \int_{\partial\Omega}
    e^{-i\kappa_\beta|\vy-\vz^s|}
    \ve_j\cdot
    \int_{\partial\Omega}
    e^{-i\kappa_\alpha|\vx-\vz^s|}
    \vu_{\mathrm{scat}}^{(j)}(\vx,\vy)
    \,d\sigma(\vx)\,d\sigma(\vy),
\end{align}
where $\alpha,\beta=P,S$. The first index $\alpha$ specifies the wave
mode used for receiver-side backpropagation in the variable $\vx$,
whereas the second index $\beta$ specifies the mode used for
source-side backpropagation in the variable $\vy$. Accordingly, the
four choices
\[
    (\alpha,\beta)
    =
    (P,P),\ (P,S),\ (S,P),\ (S,S)
\]
define the PP, PS, SP, and SS mode imaging functionals, respectively. When no confusion can arise, we simply refer to these
imaging functionals as the PP, PS, SP, and SS modes.

Our objectives are to localize the small inclusion from the measured
scattered fields
$\{\vu_{\mathrm{scat}}^{(j)}\}_{j=1,2}$, characterize the
point spread functions and resolution properties of the four imaging
functionals, and develop a procedure for identifying the orientation, the material and
geometric parameters of the inclusion.
\subsection{Elastic moment tensor}
\label{subsec:intro_EMT}

In this subsection, we briefly review the elastic moment tensor (EMT)
in two dimensions and its role in the small-volume asymptotic expansion
of the scattered field. For ease of presentation, we we restrict attention to the case of an ellipse. For general results, we refer to
\cite{ammari2007polarization,
3d546238-6f8a-3f6d-b0b4-3883234b540d}.

For an ellipse $B\subset\mathbb{R}^2$, the EMT
\[
    \sM=(m_{ijpq})_{i,j,p,q=1,2}
\]
is determined by the domain $B$, the background Lam\'e parameters
$(\lambda,\mu)$, and the inclusion Lam\'e parameters
$(\lambda',\mu')$. Its components are defined by
\begin{equation}\label{eq:EMT_definition}
    m_{ijpq}
    =
    \int_{\partial B}
    \left[
        -\frac{\partial(x_p\ve_q)}{\partial\vnu}
        +\frac{\partial(x_p\ve_q)}{\partial\vnu'}
    \right]
    \cdot\vv_{ij}\,d\sigma,
\end{equation}
where $\partial/\partial\vnu$ and $\partial/\partial\vnu'$ denote the
elastic conormal derivatives associated with $(\lambda,\mu)$ and
$(\lambda',\mu')$, respectively. The function $\vv_{ij}$ is the unique
solution of the transmission problem
\[
\renewcommand{\arraystretch}{1.5}
\left\{
\begin{aligned}
    &\fL^{\lambda,\mu}\vv_{ij}
    =0
    &&\text{in }\mathbb{R}^2\setminus\overline{B},\\
    &\fL^{\lambda',\mu'}\vv_{ij}
    =0
    &&\text{in }B,\\
    &\vv_{ij}\big|_-
    =\vv_{ij}\big|_+
    &&\text{on }\partial B,\\
    &\frac{\partial\vv_{ij}}{\partial\vnu'}\bigg|_-
    =
    \frac{\partial\vv_{ij}}{\partial\vnu}\bigg|_+
    &&\text{on }\partial B,\\
    &\vv_{ij}(\vx)-x_i\ve_j
    =\fO(|\vx|^{-1})
    &&\text{as }|\vx|\to\infty.
\end{aligned}
\right.
\]

Assume that $\rho\omega^2$ is not a Neumann eigenvalue of
$-\fL^{\lambda,\mu}$ in the inclusion. Under the small-volume
condition~\eqref{eq:small-volume-regime}, the scattered field generated
by a source at $\vx_k$ and measured at a receiver $\vx_l$ admits the
expansion
\cite{3d546238-6f8a-3f6d-b0b4-3883234b540d,
doi:10.1137/070696076}
\begin{align}\label{eq:scattered_expand}
    \vu_{\mathrm{scat}}^{(j)}(\vx_l,\vx_k)
    &:=
    \bigl[
        \vu_\varepsilon^{(j)}
        -
        \vu_{\mathrm{inc}}^{(j)}
    \bigr](\vx_l,\vx_k)
    \notag\\
    &=
    -\varepsilon^2
    \Bigl[
        \nabla_{\vz}
        \bigl(
            \mGamma^\omega(\vz,\vx_k)\ve_j
        \bigr)
        :
        \sM^\phi\nabla_{\vz}\mGamma^\omega(\vx_l,\vz)
        \notag\\
    &\hspace{3.5cm}
        +
        \omega^2(\rho-\rho')|B|
        \mGamma^\omega(\vx_l,\vz)
        \mGamma^\omega(\vz,\vx_k)\ve_j
    \Bigr]
    +
    \fO(\varepsilon^3).
\end{align}
Here, $\sM^\phi$ denotes the EMT expressed in the laboratory coordinate
system. For an inclusion whose principal axes are aligned with the
laboratory coordinate axes, $\phi=0$ and $\sM^\phi=\sM$.

The double contraction in \eqref{eq:scattered_expand} is understood
componentwise as
\begin{align}\label{eq:EMT_contraction}
    &\left[
        \nabla_{\vz}
        \bigl(
            \mGamma^\omega(\vz,\vy)\ve_j
        \bigr)
        :
        \sM^\phi\nabla_{\vz}\mGamma^\omega(\vx,\vz)
    \right]_k
    =
    \sum_{i,l,p,q=1}^{2}
    \partial_{z_i}
    \bigl(
        \mGamma^\omega(\vz,\vy)\ve_j
    \bigr)_l
    m_{ilpq}^\phi
    \partial_{z_p}
    \Gamma_{kq}^\omega(\vx,\vz).
\end{align}

The EMT possesses the following symmetries.

\begin{prop}\label{prop:EMT}
Let $\sM$ be the EMT associated with $B$, with Lam\'e parameters
$(\lambda',\mu')$ inside $B$ and $(\lambda,\mu)$ in the background.
Then
\begin{equation}\label{eq:EMT_symmetry}
    m_{ijpq}
    =
    m_{ji pq}
    =
    m_{ijqp}
    =
    m_{pqij},
    \qquad i,j,p,q=1,2.
\end{equation}
\end{prop}

Let $B$ be an
ellipse with semi-axis lengths $a$ and $b$, where $a\geq b$, and
suppose that its major and minor axes are aligned with the laboratory coordinate axes. The nonzero independent
components of its EMT in this coordinate system are given by
\begin{align}
    m_{1211}
    &=
    m_{1222}
    =
    0,\\
    m_{1111}
    &=
    |B|(\lambda+2\mu)
    \frac{
        (\mu'-\mu)(\lambda-\lambda'+\mu-\mu')
    \bigl[m^2-2(\tau-1)m\bigr]+c
    }{
       (\mu'-\mu)\bigl[\mu-\tau(\lambda'+\mu')\bigr]m^2
    +(\mu+\lambda'+\mu')(\mu+\tau\mu')
    },\\
    m_{2222}
    &=
    |B|(\lambda+2\mu)
    \frac{
        (\mu'-\mu)(\lambda-\lambda'+\mu-\mu')
    \bigl[m^2+2(\tau-1)m\bigr]+c
    }{
        (\mu'-\mu)\bigl[\mu-\tau(\lambda'+\mu')\bigr]m^2
    +(\mu+\lambda'+\mu')(\mu+\tau\mu')
    },\\
    m_{1122}
    &= |B|
\frac{
    (\lambda+2\mu)\bigl[
        (\mu'-\mu)(\lambda-\lambda'+\mu-\mu')m^2
        +(\lambda'-\lambda)(\mu'+\tau\mu)
        +(\mu'-\mu)^2
    \bigr]
}{
    (\mu'-\mu)\bigl[\mu-\tau(\lambda'+\mu')\bigr]m^2
    +(\mu+\lambda'+\mu')(\mu+\tau\mu')
},\\
    m_{1212}
    &=
   |B|
\frac{
    \mu(\mu'-\mu)(\tau+1)
}{
    (\mu'-\mu)m^2+\mu+\tau\mu'
},
\end{align}
where
\begin{align}
    c
    &=
    (\lambda'-\lambda+\mu'-\mu)(\mu+\tau\mu')
    +
    (\tau-1)(\mu'-\mu)(\mu+\lambda'+\mu'),\\
    m
    &=
    \frac{a-b}{a+b},\\
    \tau
    &=
    \frac{\lambda+3\mu}{\lambda+\mu}.
\end{align}
The remaining components are determined by the symmetry relations in
Proposition~\ref{prop:EMT}.

In the special case $m=0$, the domain $B$ is a disk, and the preceding
expressions reduce to
\begin{align}
    m_{1211}
    &=
    m_{1222}
    =
    0,
    \label{eq:m_1}\\
    m_{1111}
    &=
    m_{2222}
    =
    m_{1122}+2m_{1212},\\
    m_{1122}
    &=
    |B|
    \frac{
        (\lambda+2\mu)
        \left[
            (\lambda'-\lambda)(\mu'+\tau\mu)
            +
            (\mu'-\mu)^2
        \right]
    }{
        (\mu+\lambda'+\mu')(\mu+\tau\mu')
    },\\
    m_{1212}
    &=
    |B|
    \frac{
        \mu(\mu'-\mu)(\tau+1)
    }{
        \mu+\tau\mu'
    }.
    \label{eq:m_2}
\end{align}

When the background Lam\'e parameters are known, the four independent
EMT components in the principal-axis coordinate system can be regarded
as functions of the inclusion parameters
$(\lambda',\mu',|B|,m)$. For later use, we write
\begin{align}
    m_{1111}
    &=
    f_1(\lambda',\mu',|B|,m),
    \label{eq:f_1}\\
    m_{2222}
    &=
    f_2(\lambda',\mu',|B|,m),\\
    m_{1122}
    &=
    f_3(\lambda',\mu',|B|,m),\\
    m_{1212}
    &=
    f_4(\lambda',\mu',|B|,m).
    \label{eq:f_4}
\end{align}

We now relate the principal-axis representation of the EMT to its
representation in the laboratory coordinate system. Suppose that the
major axis of the ellipse forms an angle $\phi$ with the positive
$x_1$-axis. The corresponding transformation law is given below
\cite{ammari2007polarization}.

\begin{prop}\label{prop:rotation}
Let
\[
    \mR^\phi=(R_{ij}^\phi)_{i,j=1,2}
\]
be the rotation matrix defined in \eqref{eq:def_R_phi}. Let
\[
    \sM=(m_{klrs})_{k,l,r,s=1,2}
\]
denote the EMT expressed in the principal-axis coordinate system, and
let
\[
    \sM^\phi=(m_{ijpq}^\phi)_{i,j,p,q=1,2}
\]
denote the same tensor expressed in the laboratory coordinate system.
Then
\begin{equation}\label{eq:EMT_rotation}
    m_{ijpq}^\phi
    =
    \sum_{k,l,r,s=1}^{2}
    R_{ik}^\phi
    R_{jl}^\phi
    R_{pr}^\phi
    R_{qs}^\phi
    m_{klrs},
    \qquad i,j,p,q=1,2.
\end{equation}
\end{prop}

In Section~\ref{sec:psf_resolution}, the imaging formulas are stated
for an arbitrarily oriented ellipse and therefore explicitly involve
the orientation angle $\phi$. For clarity, their derivation is first
carried out in the principal-axis coordinate system, corresponding to
$\phi=0$. The formulas for an arbitrary orientation then follow from
the rotational transformation in Proposition~\ref{prop:rotation}.
The numerical treatment of the unknown orientation is presented in
Section~\ref{subsec:elliptic}.

\section{Point spread functions and resolution analysis}\label{sec:psf_resolution}
In this section, we derive explicit leading-order expressions for the
imaging functionals associated with the four polarization modes and
analyze their resolution properties. These results provide the
theoretical foundation for the quantitative identification procedure
developed in Section~\ref{sec:identification}.

For clarity, we first separate the contributions arising from the
density and elasticity contrasts. The corresponding small-volume
approximations of the imaging functionals are denoted by
$\fI_{\alpha,\beta}^{\mathrm{sv},D}$ and
$\fI_{\alpha,\beta}^{\mathrm{sv},E}$, respectively, where
$\alpha,\beta=P,S$; see \eqref{eq:I^D} and \eqref{eq:I^E}. When both
contrasts are present, the small-volume expansion
\eqref{eq:scattered_expand} yields
\begin{align*}
    \fI_{\alpha,\beta}(\vz^s)
    &:=
    \fI_{\alpha,\beta}^{\mathrm{D}}(\vz^s)
    +
    \fI_{\alpha,\beta}^{\mathrm{E}}(\vz^s) \notag\\
    &\simeq
    \fI_{\alpha,\beta}^{\mathrm{sv},D}(\vz^s)
    +
    \fI_{\alpha,\beta}^{\mathrm{sv},E}(\vz^s)
    =:
    \fI_{\alpha,\beta}^{\mathrm{sv}}(\vz^s),
\end{align*}

Under the scale-separation assumptions
\eqref{eq:large-frequency-distance} and
\eqref{eq:geometric-separation}, $ \fI_{\alpha,\beta}^{\mathrm{sv}} $ admits
the explicit leading-order representation
\begin{equation}\label{eq:I_asym}
    \fI_{\alpha,\beta}^{\mathrm{sv}}(\vz^s)
    \simeq
    \fI_{\alpha,\beta}^{\mathrm{asym}}(\vz^s)
    :=
    \fI_{\alpha,\beta}^{\mathrm{asym},\mathrm{D}}(\vz^s)
    +
    \fI_{\alpha,\beta}^{\mathrm{asym},\mathrm{E}}(\vz^s).
\end{equation}
Explicit expressions for
$\fI_{\alpha,\beta}^{\mathrm{asym},\mathrm{D}}$ and
$\fI_{\alpha,\beta}^{\mathrm{asym},\mathrm{E}}$ are given in
Theorems~\ref{thm:psf_density} and \ref{thm:psf_elas},
respectively.

The derivation relies on the far-field asymptotic expansions of the
fundamental solution $\mGamma^\omega$ and its derivatives summarized
in \ref{sec:farfield_asymp}. Throughout this section, $J_n$
denotes the Bessel function of the first kind of order $n$.
\subsection{Case 1: Density contrast}\label{subsec:image_den_theo}

We first consider a pure density contrast and assume that
\[
    (\lambda',\mu')=(\lambda,\mu).
\]
In this case, the leading-order scattered field is
\begin{equation}
    \bigl[
        \vu_\varepsilon^{(j)}
        -
        \vu_{\mathrm{inc}}^{(j)}
    \bigr](\vx_l,\vx_k)
    \simeq
    -\varepsilon^2\omega^2(\rho-\rho')|B|\,
    \mGamma^\omega(\vx_l,\vz)
    \mGamma^\omega(\vz,\vx_k)\ve_j.
\end{equation}
Consequently, the imaging functional can be written as, for $\alpha,\beta=P,S$,
\begin{align}\label{eq:I^D}
    \fI_{\alpha,\beta}^{\mathrm{sv},\mathrm{D}}(\vz^s)
    \simeq{}&
    -\varepsilon^2\omega^2(\rho-\rho')|B|
    \sum_{j=1}^{2}
    \int_{\partial\Omega}
    e^{-i\kappa_\beta|\vy-\vz^s|}
    \ve_j\cdot
    \int_{\partial\Omega}
    e^{-i\kappa_\alpha|\vx-\vz^s|}
    \notag\\
    &\qquad\qquad\qquad\times
    \mGamma^\omega(\vx,\vz)
    \mGamma^\omega(\vz,\vy)\ve_j
    \,d\sigma(\vx)\,d\sigma(\vy).
\end{align}

To describe the resolution properties of the imaging functional, it is useful to introduce a resolution-scaled parametrization of the search point $\vz^s$. The expected resolution of the imaging functional is governed by the Rayleigh criterion, which states that the resolution is of the order of the wavelength. Using this parametrization, we obtain explicit quantitative profiles of the imaging functionals for the different polarization modes.

\begin{thm}\label{thm:psf_density}
Assume that
\[
    (\lambda',\mu')=(\lambda,\mu),
\]
and let
\[
    \vz^s=\vz+\varepsilon^{1/2}\vr^s,
    \qquad
    r^s:=|\vr^s|,
    \qquad
    R^s:=\varepsilon^{1/2}r^s.
\]
For $\varepsilon \rightarrow 0$ and $\delta := \frac{|\vz|}{R} \rightarrow 0$, the density-contrast imaging functionals satisfy
\[
\fI_{\alpha,\beta}^{\mathrm{D}}(\vz^s) \simeq\fI_{\alpha,\beta}^{\mathrm{sv},\mathrm{D}}(\vz^s) \simeq \fI_{\alpha,\beta}^{\mathrm{asym},\mathrm{D}}(\vz^s),
\]
where
\begin{align}
    \fI_{P,P}^{\mathrm{asym},\mathrm{D}}(\vz^s)
    &=
    -\varepsilon^2
    \frac{i\pi R\omega(\rho-\rho')|B|c_P}
         {4(\lambda+2\mu)^2}
    \left(J_0^2+J_2^2\right)
    \left(\frac{\omega}{c_P}R^s\right),
    \\
    \fI_{S,S}^{\mathrm{asym},\mathrm{D}}(\vz^s)
    &=
    -\varepsilon^2
    \frac{i\pi R\omega(\rho-\rho')|B|c_S}
         {4\mu^2}
    \left(J_0^2+J_2^2\right)
    \left(\frac{\omega}{c_S}R^s\right),
    \\
    \fI_{P,S}^{\mathrm{asym},\mathrm{D}}(\vz^s)
    &=
    -\varepsilon^2
    \frac{i\pi R\omega(\rho-\rho')|B|\sqrt{c_Pc_S}}
         {4\mu(\lambda+2\mu)}
    \Bigg[
        J_0\left(\frac{\omega}{c_P}R^s\right)
        J_0\left(\frac{\omega}{c_S}R^s\right)
        -
        J_2\left(\frac{\omega}{c_P}R^s\right)
        J_2\left(\frac{\omega}{c_S}R^s\right)
    \Bigg],
    \\
    \fI_{S,P}^{\mathrm{asym},\mathrm{D}}(\vz^s)
    &=
    \fI_{P,S}^{\mathrm{asym},\mathrm{D}}(\vz^s).
\end{align}
\end{thm}

\begin{proof}
See \ref{sec:proof_I^D}.
\end{proof}

\begin{rmk}
    The results above show that, in the presence of density contrast alone, the imaging functionals are independent of the inclusion orientation angle \(\phi\) and therefore of the choice between the laboratory and principal-axis coordinate systems.
\end{rmk}

\subsection{Case 2: Elasticity contrast}
\label{subsec:image_elas_theo}

We next consider a pure elasticity contrast and assume that
\[
    \rho'=\rho.
\]
For the elliptic inclusion $D_\phi$, let $\sM^\phi$ denote the EMT expressed in the laboratory coordinate
system. The leading-order
scattered field is then
\begin{equation}
    \bigl[
        \vu_\varepsilon^{(j)}
        -
        \vu_{\mathrm{inc}}^{(j)}
    \bigr](\vx_l,\vx_k)
    \simeq
    -\varepsilon^2
    \nabla_{\vz}
    \bigl(
        \mGamma^\omega(\vz,\vx_k)\ve_j
    \bigr)
    :
    \sM^\phi
    \nabla_{\vz}\mGamma^\omega(\vx_l,\vz).
\end{equation}
Consequently, for $\alpha,\beta=P,S$,
\begin{align}\label{eq:I^E}
    \fI_{\alpha,\beta}^{\mathrm{sv},\mathrm{E}}(\vz^s)
    \simeq{}&
    -\varepsilon^2
    \sum_{j=1}^{2}
    \int_{\partial\Omega}
    e^{-i\kappa_\beta|\vy-\vz^s|}
    \ve_j\cdot
    \int_{\partial\Omega}
    e^{-i\kappa_\alpha|\vx-\vz^s|}
    \nabla_{\vz}
    \bigl(
        \mGamma^\omega(\vz,\vy)\ve_j
    \bigr)
    :
    \sM^\phi
    \nabla_{\vz}\mGamma^\omega(\vx,\vz)
    \,d\sigma(\vx)\,d\sigma(\vy).
\end{align}

Unlike the contribution arising from density contrast, the contribution due to elasticity contrast depends on the orientation of the inclusion. It is therefore natural to express the resulting imaging profiles in the principal-axis coordinate system of the inclusion. Using the same resolution-scaled parametrization as in Section~\ref{subsec:image_den_theo}, we derive below the quantitative profiles of the imaging functionals for the four polarization modes in this coordinate system.

\begin{thm}\label{thm:psf_elas}
Assume that $\rho'=\rho$, and the elliptic inclusion $D_\phi$ has orientation angle $\phi$. Let
\[
    \vz^s
    =
    \vz+\varepsilon^{1/2}\vr^s,
    \qquad
    \vr^s=(r_1^s,r_2^s),
\]
where $r_j^s = \ve^\top_j\vr^s$, for $j=1,2$,
and let
\[
    r^s:=|\vr^s|,
    \quad
    R^s:=\varepsilon^{1/2}r^s=|\vz^s-\vz|,
    \quad
    \varphi^s
    :=
    \arg(r_1^s+i r_2^s).
\]
We then define the relative sampling angle
\begin{equation}\label{eq:relative_sampling_angle}
    \theta^s
    :=
    \varphi^s-\phi
    \quad (\operatorname{mod}2\pi),
\end{equation}
For $\varepsilon \rightarrow 0$ and $\delta := \frac{|\vz|}{R} \rightarrow 0$, the elasticity-contrast imaging functionals satisfy
\[
\fI_{\alpha,\beta}^{\mathrm{E}}(\vz^s) \simeq\fI_{\alpha,\beta}^{\mathrm{sv},\mathrm{E}}(\vz^s) \simeq \fI_{\alpha,\beta}^{\mathrm{asym},\mathrm{E}}(\vz^s).
\]
Expressed in the principal-axis coordinate system in terms of the
relative sampling angle $\theta^s$, the explicit asymptotic profiles
are given by
\begin{align}
    \fI_{P,P}^{\mathrm{asym},\mathrm{E}}(\vz^s)
    =&
    -\frac{i\kappa_P\pi R\varepsilon^2}
    {8(\lambda+2\mu)^2}
    \Bigg[
        A_{P,P}^{(1)}(\theta^s)
        J_1^2\left(\frac{\omega}{c_P}R^s\right)
        +
        A_{P,P}^{(2)}(\theta^s)
        J_3^2\left(\frac{\omega}{c_P}R^s\right)
        +
        A_{P,P}^{(3)}(\theta^s)
        J_1\left(\frac{\omega}{c_P}R^s\right)
        J_3\left(\frac{\omega}{c_P}R^s\right)
    \Bigg],
    \label{eq:PP_elasticity_profile}
    \\
    \fI_{S,S}^{\mathrm{asym},\mathrm{E}}(\vz^s)
    =&
    -\frac{i\kappa_S\pi R\varepsilon^2}
    {8\mu^2}
    \Bigg[
        A_{S,S}^{(1)}(\theta^s)
        J_1^2\left(\frac{\omega}{c_S}R^s\right)
        +
        A_{S,S}^{(2)}(\theta^s)
        J_3^2\left(\frac{\omega}{c_S}R^s\right)
        +
        A_{S,S}^{(3)}(\theta^s)
        J_1\left(\frac{\omega}{c_S}R^s\right)
        J_3\left(\frac{\omega}{c_S}R^s\right)
    \Bigg],
    \label{eq:SS_elasticity_profile}
    \\
    \fI_{P,S}^{\mathrm{asym},\mathrm{E}}(\vz^s)
    =&
    -\frac{
        i\sqrt{\kappa_S\kappa_P}\,\pi R\varepsilon^2
    }{
        8\mu(\lambda+2\mu)
    }
    \Bigg[
        A_{P,S}^{(1)}(\theta^s)
        J_1\left(\frac{\omega}{c_P}R^s\right)
        J_1\left(\frac{\omega}{c_S}R^s\right)
        +
        A_{P,S}^{(2)}(\theta^s)
        J_3\left(\frac{\omega}{c_P}R^s\right)
        J_3\left(\frac{\omega}{c_S}R^s\right)
        \notag\\
    &\qquad
        +
        A_{P,S}^{(3)}(\theta^s)
        J_1\left(\frac{\omega}{c_S}R^s\right)
        J_3\left(\frac{\omega}{c_P}R^s\right)
        +
        A_{P,S}^{(4)}(\theta^s)
        J_1\left(\frac{\omega}{c_P}R^s\right)
        J_3\left(\frac{\omega}{c_S}R^s\right)
    \Bigg],
    \label{eq:PS_elasticity_profile}
    \\
    \fI_{S,P}^{\mathrm{asym},\mathrm{E}}(\vz^s)
    =&\fI_{P,S}^{\mathrm{asym},\mathrm{E}}(\vz^s).
    \label{eq:SP_elasticity_profile}
\end{align}
The angular coefficients in the PP mode are
\begin{align}
    A_{P,P}^{(1)}(\theta^s)
    &=
    M^{(1)}
    +
    M^{(2)}\cos2\theta^s,
    \label{eq:A_PP_1}\\
    A_{P,P}^{(2)}(\theta^s)
    &=
    M^{(3)},
    \label{eq:A_PP_2}\\
    A_{P,P}^{(3)}(\theta^s)
    &=
    M^{(4)}\cos4\theta^s
    -
    M^{(2)}\cos2\theta^s.
    \label{eq:A_PP_3}
\end{align}
The angular coefficients in the SS mode are
\begin{align}
    A_{S,S}^{(1)}(\theta^s)
    &=
    A_{S,S}^{(2)}(\theta^s)
    =
    M^{(3)},
    \label{eq:A_SS_12}\\
    A_{S,S}^{(3)}(\theta^s)
    &=
    -M^{(4)}\cos4\theta^s.
    \label{eq:A_SS_3}
\end{align}
The angular coefficients in the SP mode are
\begin{align}
    A_{S,P}^{(1)}(\theta^s)
    &=
    M^{(3)}
    +
    \frac{M^{(2)}}{2}\cos2\theta^s,
    \label{eq:A_SP_1}\\
    A_{S,P}^{(2)}(\theta^s)
    &=
    -M^{(3)},
    \label{eq:A_SP_2}\\
    A_{S,P}^{(3)}(\theta^s)
    &=
    \frac{M^{(4)}}{2}\cos4\theta^s,
    \label{eq:A_SP_3}\\
    A_{S,P}^{(4)}(\theta^s)
    &=
    -\frac{M^{(4)}}{2}\cos4\theta^s
    +
    \frac{M^{(2)}}{2}\cos2\theta^s.
    \label{eq:A_SP_4}
\end{align}
The corresponding PS coefficients satisfy
\begin{equation}\label{eq:A_PS_SP_relation}
    A_{P,S}^{(j)}(\theta^s)
    =
    A_{S,P}^{(j)}(\theta^s),
    \qquad j=1,\ldots,4.
\end{equation}

Here, $m_{1111}$, $m_{2222}$, $m_{1122}$, and $m_{1212}$ denote the
components of the EMT expressed in the principal-axis coordinate
system of the ellipse, and
\begin{align}
    M^{(1)}
    &=
    m_{1212}
    +
    \frac{3}{2}m_{1122}
    +
    \frac{5}{4}
    \bigl(m_{1111}+m_{2222}\bigr),
    \label{eq:M_1}\\
    M^{(2)}
    &=
    m_{1111}-m_{2222},
    \label{eq:M_2}\\
    M^{(3)}
    &=
    \frac{1}{4}
    \bigl(m_{1111}+m_{2222}\bigr)
    +
    m_{1212}
    -
    \frac{1}{2}m_{1122},
    \label{eq:M_3}\\
    M^{(4)}
    &=
    2m_{1212}
    +
    m_{1122}
    -
    \frac{1}{2}
    \bigl(m_{1111}+m_{2222}\bigr).
    \label{eq:M_4}
\end{align}

In the special case where $B$ is a disk,
\begin{equation}\label{eq:M_disk}
    M^{(1)}
    =
    6m_{1212}+4m_{1122},
    \qquad
    M^{(3)}
    =
    2m_{1212},
    \qquad
    M^{(2)}
    =
    M^{(4)}
    =
    0.
\end{equation}
Consequently, the angular dependence disappears, consistently with the
rotational symmetry of the disk.
\end{thm}

\begin{proof}
The formulas are first derived for $\phi=0$, in which case the
principal axes of the ellipse are aligned with the laboratory
coordinate axes; see \ref{sec:proof_I^E}. For an arbitrary
orientation angle $\phi$, by introducing the relative sampling angle
\[
    \theta^s=\varphi^s-\phi,
\]
the formulas stated in the theorem are derived.
\end{proof}

\begin{rmk}\label{rmk:coordinate_representations}
The above results show that the elasticity-contrast imaging
functionals depend on the orientation of the inclusion. In
Theorem~\ref{thm:psf_elas}, they are expressed in the principal-axis
coordinate system in terms of $(\theta^s,\sM)$. Equivalently, they may
be expressed in the laboratory coordinate system in terms of
$(\varphi^s,\sM^\phi)$, where $\sM^\phi$ is obtained from $\sM$ through
the rotation formula in Proposition~\ref{prop:rotation}.
\end{rmk}

For imaging purposes, all four imaging
functionals provide an estimate of the inclusion location. We define
the resolution of an imaging functional in terms of the characteristic
width of its point-spread function. The asymptotic profiles established
in the preceding theorems show that the inclusion can be localized by
applying Kirchhoff migration to the scattered-field data and further
reveal the resolution properties of the different polarization modes.
The resolution of the PP and SS modes is governed primarily by the
\(P\)- and \(S\)-wave wavelengths, respectively, whereas the mixed PS
and SP modes depend on the interaction between the two wavelengths.
Since \(c_S<c_P\), the \(S\)-wave wavelength is shorter than the
\(P\)-wave wavelength at a fixed excitation frequency. Consequently,
the SS mode has a narrower characteristic resolution scale than the PP
mode and therefore provides higher spatial resolution.

The imaging functionals
\(\{\fI_{\alpha,\beta}^{\mathrm{asym},\mathrm{E}}\}_{\alpha,\beta=P,S}\) are expressed as
linear combinations of products of Bessel functions, with coefficients
determined by the EMT of the inclusion. In
particular, the PP mode  contains four independent EMT
coefficients, \(\{M^{(j)}\}_{j=1}^4\), whereas each of the SS, PS, and
SP modes contains only two independent combinations of
the EMT components. This distinction is essential for the
identification problem. Among the four polarization modes, only the PP
mode contains sufficient independent information to uniquely recover
the EMT and, consequently, the material and geometric parameters of an
elliptic inclusion. The remaining modes are useful for localization
but do not, individually, provide sufficient information for
quantitative identification. The resulting identification procedure is
developed in Section~\ref{sec:identification}.

\section{On the identification of inclusion properties}\label{sec:identification}

In this section, we investigate the inverse problem of identifying the properties of the inclusion, including its Lam\'e parameters, density, and geometric parameters. Section~\ref{subsec:elliptic} is devoted to the case of an elliptic inclusion, while Section~\ref{subsec:spherical} considers the special case of a disk-shaped inclusion.

\subsection{Elliptic inclusion}\label{subsec:elliptic}

In this section, we first establish a uniqueness result for the identification of the EMT from the PP mode. Then, we present the practical implementation of the identification procedure.

The asymptotic analysis of the previous section shows that among the four polarization modes, only the PP mode contains sufficient independent information for recovering all inclusion parameters. Recall that, in the presence of both density and elasticity contrasts, adopting the parametrization in Theorem~\ref{thm:psf_density} and Theorem~\ref{thm:psf_elas}, the PP mode admits the asymptotic representation in terms of principal-axis coordinates:
\begin{align}
\label{eq:both_contrast}
\fI_{P,P}(\vz^s) \simeq \fI^{\mathrm{asym}}_{P,P}(\vz^s)=
&-\varepsilon^2\frac{i\pi R\omega M^{(0)} c_P}{4(\lambda+2\mu)^2}
   \bra{J_0^2+J_2^2}
   \bra{\frac{\omega}{c_P}R^s}
   \notag\\
&\quad-\frac{i\kappa_P\pi R\varepsilon^2}{8(\lambda+2\mu)^2}
   \Bigl[
      (M^{(1)}+M^{(2)}\cos 2\theta^s)
      J_1^2\!\Bigl(\frac{\omega}{c_P}R^s\Bigr)
      +M^{(3)}
      J_3^2\!\Bigl(\frac{\omega}{c_P}R^s\Bigr)
      \notag\\
&\quad
      +(M^{(4)}\cos4\theta^s-M^{(2)}\cos2\theta^s)
      J_1\!\Bigl(\frac{\omega}{c_P}R^s\Bigr)
      J_3\!\Bigl(\frac{\omega}{c_P}R^s\Bigr)
   \Bigr],
\end{align}
where $M^{(0)}=(\rho-\rho')|B|$ and $\{M^{(j)}\}_{j=1}^4$ are defined in Theorem~\ref{thm:psf_elas}.

The following theorem establishes the unique determination of the EMT from the PP mode.

\begin{thm}\label{thm:emt_uniqueness}
Assume that the center $\vz$ and the orientation angle $\phi$ of the
elliptic inclusion have been identified. Consider sampling points of the form
\[
    \vz_{m,n}^s
    =
    \vz
    +
    R_m
    \begin{pmatrix}
        \cos\varphi_n^s\\
        \sin\varphi_n^s
    \end{pmatrix},
    \qquad
    m=1,\ldots,N_R,
    \quad
    n=1,\ldots,N_{\varphi^s},
\]
where $\varphi_n^s$ denotes the sampling angle in the laboratory
coordinate system. Define the corresponding angle in the
principal-axis coordinate system by
\[
    \theta_n^s
    :=
    \varphi_n^s-\phi
    \quad(\operatorname{mod}2\pi).
\]
Suppose that the angular sampling matrix
\[
    \mV
    =
    \begin{pmatrix}
        1 & \cos(2\theta_1^s) & \cos(4\theta_1^s)\\
        \vdots & \vdots & \vdots\\
        1 & \cos(2\theta_{N_{\varphi^s}}^s)
          & \cos(4\theta_{N_{\varphi^s}}^s)
    \end{pmatrix}.
\]
has rank three. Define
\[
\begin{aligned}
    g_0(R)
    &:=
    \bigl[J_0^2+J_2^2\bigr]
    \left(\frac{\omega}{c_P}R\right),\\
    g_1(R)
    &:=
    J_1^2\left(\frac{\omega}{c_P}R\right),\\
    g_3(R)
    &:=
    J_3^2\left(\frac{\omega}{c_P}R\right),
\end{aligned}
\]
and assume that three sampling radii, relabeled as $R_1,R_2,R_3$,
satisfy
\[
    \det
    \begin{pmatrix}
        g_0(R_1) & g_1(R_1) & g_3(R_1)\\
        g_0(R_2) & g_1(R_2) & g_3(R_2)\\
        g_0(R_3) & g_1(R_3) & g_3(R_3)
    \end{pmatrix}
    \neq0.
\]
Suppose, in addition, that the set of sampling radii contains a radius
$R_{m_2}$ such that
\[
    J_1^2\left(\frac{\omega}{c_P}R_{m_2}\right)
    -
    J_1\left(\frac{\omega}{c_P}R_{m_2}\right)
    J_3\left(\frac{\omega}{c_P}R_{m_2}\right)
    \neq0,
\]
and a radius $R_{m_4}$ such that
\[
    J_1\left(\frac{\omega}{c_P}R_{m_4}\right)
    J_3\left(\frac{\omega}{c_P}R_{m_4}\right)
    \neq0.
\]
Then the coefficients
$M^{(0)},M^{(1)},\ldots,M^{(4)}$, and hence the EMT expressed in the
principal-axis coordinate system, are uniquely determined from the
leading-order PP imaging functional sampled at
$\{\vz_{m,n}^s\}$.
\end{thm}

\begin{proof}
See \ref{proof:emt_uniqueness}.
\end{proof}

\begin{rmk}\label{rem:sampling_conditions}
The sampling conditions in Theorem~\ref{thm:emt_uniqueness} are not
restrictive in practice. A uniform angular grid with
$N_{\theta^s}\geq9$ provides a convenient sufficient condition for
avoiding aliasing of the Fourier modes up to order four. The excluded
radial configurations correspond only to isolated zeros or exceptional
choices of the sampling radii. In the ideal case, the proof indicates that
five suitably chosen sampling points are sufficient: three points on
one sampling circle separate the angular modes, and two additional
points along one sampling direction provide the remaining radial
information. In practice, an overdetermined set of radial and angular
samples is used, and the EMT coefficients are recovered by least
squares, which ensures the robustness of the estimate in the presence of model errors or measurement errors.
\end{rmk}

Once the EMT has been recovered in the principal-axis coordinate
system of the ellipse, in which the coordinate axes are aligned with the major and
minor axes of the ellipse, the Lam\'e parameters and geometric parameters
$(\lambda',\mu',|B|,m)$ can be obtained by solving the nonlinear system
\eqref{eq:f_1}--\eqref{eq:f_4}. The density is subsequently recovered
from
\[
    \rho'
    =
    \rho-\frac{M^{(0)}}{|B|}.
\]
From a practical perspective, the inclusion area $|B|$ and ellipticity
parameter $m$ are often more relevant than the individual semi-axes.
The semi-axes can be uniquely reconstructed from these quantities
according to
\[
    a
    =
    \sqrt{\frac{(1+m)|B|}{(1-m)\pi}},
    \qquad
    b
    =
    \sqrt{\frac{(1-m)|B|}{(1+m)\pi}}.
\]
Accordingly, the reconstruction procedure focuses on recovering
$(|B|,m)$ rather than $(a,b)$ directly.

As indicated by Theorem~\ref{thm:emt_uniqueness}, recovering the EMT and the inclusion parameters requires prior estimates of the inclusion center \(\vz\) and orientation angle \(\phi\). We now describe the practical identification procedure, beginning with the localization of the inclusion center from the measured scattered field.

Let $\fI^{\mathrm{meas}}_{\alpha,\beta}$ denote the imaging functional with
polarization modes $(\alpha,\beta)$ computed from the measured
scattered field
\[
    \vu_{\mathrm{scat}}^{(j)}
    =
    \vu_\varepsilon^{(j)}-\vu_{\mathrm{inc}}^{(j)}
\]
on a sampling grid $\vLambda$. Since all four polarization modes
provide an estimate of the inclusion location, any one of them can, in
principle, be used for localization. We define
\[
    \widehat{\vz}
    =
    \operatorname*{arg\,max}_{\vz^s\in\vLambda}
    \left|
        \operatorname{Im}
        \bigl(\fI^{\mathrm{meas}}_{S,S}(\vz^s)\bigr)
    \right|,
    \qquad \alpha,\beta=P,S.
\]
The imaginary part is used because the leading-order asymptotic
formulas predict purely imaginary imaging functionals. As shown in
Section~\ref{sec:psf_resolution}, the SS mode has the highest spatial
resolution among the four modes and is therefore preferred for
localizing the inclusion center. For simplicity, we assume that the imaging functional has a dominant peak at the inclusion center, as in configurations dominated by the density contribution. Localization in more general cases is discussed in Remark~\ref{rmk:general_case}.

The preceding formulas were derived in the principal-axis coordinate
system. For an elliptic inclusion with unknown orientation, the principal axes are not aligned with the laboratory coordinate axes. The orientation angle $\phi$ is therefore unknown and must be estimated from the imaging functionals.

To this end, let $\varphi\in[0,2\pi)$ denote the polar angle measured from the
$x_1$-axis. At a fixed sampling radius $R_*>0$, we define the function 
\[
\fP(\varphi):= \operatorname{Im}\left[\fI_{P,P}^{\mathrm{meas}}(\hat{\vz}+R_*(\cos\varphi,\sin\varphi))\right].
\]
By the asymptotic profile \eqref{eq:both_contrast}, we have that
$\mathcal{P}$ takes the form
\begin{equation}\label{eq:orientation_shifted}
    \mathcal{P}(\varphi)
    \simeq
    a_0
    +a_2\cos2(\varphi-\phi)
    +a_4\cos4(\varphi-\phi),
\end{equation}
where $a_0$, $a_2$, and $a_4$ are real-valued amplitudes depending on
$R_*$ and the EMT coefficients.

Using
\[
    \cos2(\varphi-\phi)
    =
    \cos2\varphi\cos2\phi
    +
    \sin2\varphi\sin2\phi
\]
and
\[
    \cos4(\varphi-\phi)
    =
    \cos4\varphi\cos4\phi
    +
    \sin4\varphi\sin4\phi,
\]
we rewrite \eqref{eq:orientation_shifted} as
\begin{align}\label{eq:orientation_fourier}
    \mathcal{P}(\varphi)
    \simeq{}&
    a_0
    +a_{2,c}\cos2\varphi
    +a_{2,s}\sin2\varphi
    \notag\\
    &+
    a_{4,c}\cos4\varphi
    +a_{4,s}\sin4\varphi,
\end{align}
where
\[
    a_{2,c}:=a_2\cos2\phi,
    \qquad
    a_{2,s}:=a_2\sin2\phi,
\]
and
\[
    a_{4,c}:=a_4\cos4\phi,
    \qquad
    a_{4,s}:=a_4\sin4\phi.
\]

The five Fourier coefficients
\[
    a_0,\quad
    a_{2,c},\quad
    a_{2,s},\quad
    a_{4,c},\quad
    a_{4,s}
\]
are estimated from uniformly distributed angular samples of
$\mathcal{P}$ by linear least squares. Their estimates are denoted by
\[
    \widehat{a}_0,\quad
    \widehat{a}_{2,c},\quad
    \widehat{a}_{2,s},\quad
    \widehat{a}_{4,c},\quad
    \widehat{a}_{4,s}.
\]

Provided that
\[
    \hat{a}_{2,c}^2+\hat{a}_{2,s}^2\neq0,
\]
the second-order angular harmonic yields the preliminary orientation
estimate
\begin{equation}\label{eq:orientation_estimate}
    \widehat{\phi}_0
    :=
    \frac{1}{2}
    \arg\!\left(
        \widehat{a}_{2,c}
        +i\widehat{a}_{2,s}
    \right).
\end{equation}
Since an ellipse is invariant under a rotation by $\pi$, its
orientation is intrinsically defined modulo $\pi$. Moreover, because
the sign of the second-order harmonic amplitude is not known a priori,
\eqref{eq:orientation_estimate} gives two candidate orientations,
\[
    \widehat{\phi}_0
    \quad\text{and}\quad
    \widehat{\phi}_0+\frac{\pi}{2}
    \qquad (\operatorname{mod}\pi),
\]
corresponding to an interchange of the two principal axes. The ambiguity of principal-axis direction can be resolved in the
subsequent parameter-identification step by selecting the candidate
for which the first principal direction corresponds to the recovered
major semi-axis, under the convention $a\geq b$.

The fourth-order coefficients satisfy
\[
    a_{4,c}=a_4\cos4\phi,
    \qquad
    a_{4,s}=a_4\sin4\phi,
\]
and their estimates may be used to check the consistency of the
orientation obtained from \eqref{eq:orientation_estimate}.

To proceed, we introduce the rotation matrix
\[
    \mR^{\widehat{\phi}}
    =
    \begin{pmatrix}
        \cos\widehat{\phi} & -\sin\widehat{\phi}\\
        \sin\widehat{\phi} &  \cos\widehat{\phi}
    \end{pmatrix}.
\]
For each selected sampling point $\vz_k^s\in\vLambda$, we transform
its position relative to the estimated center into the estimated
principal-axis coordinate system:
\[
    \vy_k^s
    =
    \bigl(\mR^{\widehat{\phi}}\bigr)^\top
    (\vz_k^s-\widehat{\vz}).
\]
We then define the corresponding polar coordinates and scaled radial
variable by
\[
    R_k^s:=|\vy_k^s|,
    \qquad
    \theta_k^s:=\arg(y_{k,1}^s+i y_{k,2}^s),
    \qquad
    t_k:=\frac{\omega}{c_P}R_k^s.
\]
Thus, $\theta_k^s$ is the polar angle of the sampling direction
measured relative to the first estimated principal axis. Since the
orientation estimated from the second-order angular harmonic is
initially determined only modulo $\pi/2$, the two principal axes may be
interchanged. This ambiguity is resolved after parameter reconstruction
by ordering the recovered semi-axes according to $a\geq b$ and
assigning $\widehat{\phi}$ to the major-axis direction.

For convenience, define
\[
    C_0
    :=
    -\frac{i\pi R\omega c_P\varepsilon^2}
    {4(\lambda+2\mu)^2},
    \qquad
    C_1
    :=
    -\frac{i\kappa_P\pi R\varepsilon^2}
    {8(\lambda+2\mu)^2}.
\]
By collecting the coefficients of
$M^{(0)},\ldots,M^{(4)}$ in \eqref{eq:both_contrast}, the
leading-order PP mode at $\vz_k^s$ can be expressed
as
\begin{equation}\label{eq:linear_imaging_profile}
    \fI_{P,P}^{\mathrm{asym}}(\vz_k^s)
    =
    \vd_k^\top\vm,
\end{equation}
where
\[
    \vm
    :=
    \begin{pmatrix}
        M^{(0)}\\
        M^{(1)}\\
        M^{(2)}\\
        M^{(3)}\\
        M^{(4)}
    \end{pmatrix},
\]
and the column vector $\vd_k\in\mathbb{C}^5$ is defined by
\begin{equation}\label{eq:design_vector}
    \vd_k
    :=
    \begin{pmatrix}
        C_0\bigl[J_0^2(t_k)+J_2^2(t_k)\bigr]
        \\[2mm]
        C_1J_1^2(t_k)
        \\[2mm]
        C_1\cos(2\theta_k^s)
        \bigl[J_1^2(t_k)-J_1(t_k)J_3(t_k)\bigr]
        \\[2mm]
        C_1J_3^2(t_k)
        \\[2mm]
        C_1\cos(4\theta_k^s)J_1(t_k)J_3(t_k)
    \end{pmatrix}.
\end{equation}

For $K$ selected sampling points, we define the design matrix and the
observation vector by
\begin{equation}\label{eq:design_matrix}
    \mM_{\mathrm{design}}
    :=
    \begin{pmatrix}
        \vd_1^\top\\
        \vd_2^\top\\
        \vdots\\
        \vd_K^\top
    \end{pmatrix}
    \in\mathbb{C}^{K\times5},
    \qquad
    \vb
    :=
    \begin{pmatrix}
        \fI_{P,P}^{\mathrm{meas}}(\vz_1^s)\\
        \fI_{P,P}^{\mathrm{meas}}(\vz_2^s)\\
        \vdots\\
        \fI_{P,P}^{\mathrm{meas}}(\vz_K^s)
    \end{pmatrix}
    \in\mathbb{C}^{K}.
\end{equation}
The imaging coefficients are estimated by solving
\begin{equation}\label{eq:coefficient_least_squares}
    \widehat{\vm}
    =
    \operatorname*{arg\,min}_{\vm\in\mathbb{R}^5}
    \left\|
        \mM_{\mathrm{design}}\vm-\vb
    \right\|_2^2,
    \qquad
    \widehat{\vm}
    =
    \begin{pmatrix}
        \widehat M^{(0)}\\
        \widehat M^{(1)}\\
        \widehat M^{(2)}\\
        \widehat M^{(3)}\\
        \widehat M^{(4)}
    \end{pmatrix}.
\end{equation}
Although only five coefficients are unknown, using an overdetermined
set of sampling points improves the stability of the reconstruction.

Since the angular coordinates used to construct
$\mM_{\mathrm{design}}$ are measured in the estimated principal-axis
coordinate system, the recovered coefficients
$\{\widehat M^{(j)}\}_{j=1}^4$ determine the EMT components in the same
coordinate system. Specifically,
\[
    \widehat{\vv}
    =
    \mM_c^{-1}
    \begin{pmatrix}
        \widehat M^{(1)}\\
        \widehat M^{(2)}\\
        \widehat M^{(3)}\\
        \widehat M^{(4)}
    \end{pmatrix},
    \qquad
    \widehat{\vv}
    =
    \begin{pmatrix}
        \widehat m_{1111}\\
        \widehat m_{2222}\\
        \widehat m_{1122}\\
        \widehat m_{1212}
    \end{pmatrix},
\]
where $\mM_c$ is defined in \eqref{eq:coeff_to_EMT}. Since
\[
    \operatorname{cond}(\mM_c)\approx2.90,
\]
this linear inversion is numerically well conditioned. No further
rotation of the recovered EMT is required before identifying the
material and geometric parameters, because the relations
\eqref{eq:f_1}--\eqref{eq:f_4} are formulated in the principal-axis
coordinate system.

The unknown parameters $(\lambda',\mu',|B|,m)$ are subsequently
recovered by nonlinear least squares. Since these parameters may differ
by several orders of magnitude, we introduce the logarithmic
reparameterization
\[
    \lambda'_{\log}=\log\lambda',
    \qquad
    \mu'_{\log}=\log\mu',
    \qquad
    |B|_{\log}=\log|B|.
\]
For a noncircular ellipse with $a>b>0$, the ellipticity parameter
satisfies $0<m<1$. To enforce this constraint, we introduce an
unconstrained variable $m_{\mathrm{logi}}\in\mathbb R$ through
\[
    m_{\mathrm{logi}}
    =
    \log\left(\frac{m}{1-m}\right),
    \qquad
    m=\frac{1}{1+e^{-m_{\mathrm{logi}}}}.
\]
Together with the logarithmic parametrizations of $\lambda'$, $\mu'$,
and $|B|$, this yields an unconstrained optimization problem.

Define
\begin{align}\label{eq:def_F}
\mF\left(
\lambda'_{\log},
\mu'_{\log},
|B|_{\log},
m_{\mathrm{logi}}
\right)
:=
\begin{pmatrix}
f_1\left(
e^{\lambda'_{\log}},
e^{\mu'_{\log}},
e^{|B|_{\log}},
\frac{1}{1+e^{-m_{\mathrm{logi}}}}
\right)
\\
f_2\left(
e^{\lambda'_{\log}},
e^{\mu'_{\log}},
e^{|B|_{\log}},
\frac{1}{1+e^{-m_{\mathrm{logi}}}}
\right)
\\
f_3\left(
e^{\lambda'_{\log}},
e^{\mu'_{\log}},
e^{|B|_{\log}},
\frac{1}{1+e^{-m_{\mathrm{logi}}}}
\right)
\\
f_4\left(
e^{\lambda'_{\log}},
e^{\mu'_{\log}},
e^{|B|_{\log}},
\frac{1}{1+e^{-m_{\mathrm{logi}}}}
\right)
\end{pmatrix},
\end{align}
where $\{f_j\}_{j=1}^4$ are defined in
\eqref{eq:f_1}--\eqref{eq:f_4}. The parameter reconstruction is
formulated as
\begin{align}\label{eq:nls}
&\left(
\log\widehat{\lambda}',
\log\widehat{\mu}',
\log\widehat{|B|},
\widehat m_{\mathrm{logi}}
\right)
 =
\operatorname*{arg\,min}_{
(\lambda'_{\log},\mu'_{\log},|B|_{\log},m_{\mathrm{logi}})
\in\mathbb R^4
}
\left\|
    \widehat{\vv}
    -
    \mF\left(
      \lambda'_{\log},
      \mu'_{\log},
      |B|_{\log},
      m_{\mathrm{logi}}
    \right)
\right\|_2^2.
\end{align}
The ellipticity parameter is then recovered as
\[
    \widehat m
    =
    \frac{1}{1+e^{-\widehat m_{\mathrm{logi}}}}.
\]
Finally, the inclusion density is obtained from
\[
    \widehat{\rho}'
    =
    \rho-\frac{\widehat M^{(0)}}{\widehat{|B|}},
\]
and the semi-axes are reconstructed according to
\[
    \widehat a
    =
    \sqrt{
      \frac{(1+\widehat m)\widehat{|B|}}
           {(1-\widehat m)\pi}
    },
    \qquad
    \widehat b
    =
    \sqrt{
      \frac{(1-\widehat m)\widehat{|B|}}
           {(1+\widehat m)\pi}
    }.
\]
The resulting unconstrained least-squares problem is solved using
a multistart quasi-Newton method. The complete reconstruction
procedure is summarized in Algorithm~\ref{alg:inclusion}.

\begin{rmk}\label{rmk:general_case}
When the density and elasticity contributions interfere, the largest
magnitude of the imaging functional need not occur at the inclusion
center. Nevertheless, the leading-order SS profiles remain
centrally symmetric about $\vz$, since their angular dependence
contains only even harmonics. This observation motivates locating the center by
matching opposite parts of the image.

For the ideal leading-order profile, define
\[
\mathcal E(\vx)
:=
\left|\operatorname{Im}
\fI_{S,S}^{\mathrm{asym}}(\vx)\right|^2,
\qquad
\mathcal C(\boldsymbol s)
:=
\int_{\mathbb R^2}
\mathcal E(\vx)\mathcal E(\boldsymbol s-\vx)\,d\vx.
\]
The quantity $\mathcal C(2\boldsymbol c)$ measures the overlap of the
image energy with its reflection about a candidate center
$\boldsymbol c$. Since $\mathcal E(2\vz-\vx)=\mathcal E(\vx)$, the
Cauchy--Schwarz inequality yields, for $\mathcal E\in L^2(\mathbb R^2)$,
\[
\mathcal C(\boldsymbol s)
\leq
\|\mathcal E\|_{L^2}^2
=
\mathcal C(2\vz).
\]
This motivates estimating the center as half the location of a
maximum of $\mathcal C$. Applied to $\operatorname{Im}(\fI_{S,S}^{\mathrm{meas}})$, the same
construction provides an initial center estimate. This estimate may
be further refined by fitting
$\operatorname{Im}(\fI_{S,S}^{\mathrm{meas}})$ and $\operatorname{Im}(\fI_{S,S}^{\mathrm{asym}})$ numerically, with the center treated
as an unknown translation.
\end{rmk}

\begin{algorithm}[p]
\small
\caption{Identification of the properties of an elliptic inclusion
from scattered-field data}
\label{alg:inclusion}

\begin{algorithmic}[1]

\STATE \textbf{Input:}
Scattered-field data
$\{\vu_{\mathrm{scat}}^{(j)}\}_{j=1,2}$,
sampling grid $\vLambda$, and known background parameters
$(\lambda,\mu,\rho)$.

\STATE \textbf{Output:}
Estimated properties
$(\widehat{\vz},\widehat{\phi},
\widehat{\lambda}',\widehat{\mu}',
\widehat{|B|},\widehat m,\widehat{\rho}')$.

\STATE Compute $\fI^{\mathrm{meas}}_{P,P}$ and $\fI^{\mathrm{meas}}_{S,S}$ on
$\vLambda$ using \eqref{eq:imaging_func}.

\STATE Estimate the inclusion center:
\[
    \widehat{\vz}
    \gets
    \operatorname*{arg\,max}_{\vx_k\in\vLambda}
    \left|
      \operatorname{Im}
      \bigl(\fI^{\mathrm{meas}}_{S,S}(\vx_k)\bigr)
    \right|.
\]

\STATE Choose a radius $R_*$ and evaluate the angular profile
\[
    \mathcal P(\varphi)
    \gets
    \operatorname{Im}
    \left[
      \fI^{\mathrm{meas}}_{P,P}
      \bigl(
        \widehat{\vz}+R_*(\cos\varphi,\sin\varphi)
      \bigr)
    \right].
\]

\STATE Fit \eqref{eq:orientation_fourier} by linear least squares and
estimate the orientation:
\[
    \widehat{\phi}
    \gets
    \frac12
    \operatorname{atan2}
    \bigl(\widehat a_{2,s},\widehat a_{2,c}\bigr)
    \quad
    \left(\operatorname{mod}\frac{\pi}{2}\right),
\]
where $\operatorname{atan2}(y,x)$ denotes the polar angle of the vector $(x,y)$.

\FORALL{$\vx_k\in\vLambda$ or a selected subset}
    \STATE Set
    $\vy_k^s\gets
    (\mR^{\widehat{\phi}})^\top
    (\vx_k-\widehat{\vz})$.
    \STATE Compute
    $R_k^s\gets|\vy_k^s|$ and
    $\theta_k^s\gets\arg(y_{k,1}^s+i y_{k,2}^s)$.
\ENDFOR

\STATE Evaluate the Bessel-function combinations in
\eqref{eq:both_contrast} at $(R_k^s,\theta_k^s)$ and construct
$\mM_{\mathrm{design}}$.

\STATE Assemble
$\vb\gets
(\fI^{\mathrm{meas}}_{P,P}(\vx_1),\ldots,
 \fI^{\mathrm{meas}}_{P,P}(\vx_K))^\top$.

\STATE Recover the imaging coefficients:
\[
    \widehat{\vm}
    \gets
    \mM_{\mathrm{design}}^\dagger\vb,
    \qquad
    \widehat{\vm}
    =
    \bigl(
      \widehat M^{(0)},\widehat M^{(1)},\ldots,\widehat M^{(4)}
    \bigr)^\top,
\]
where $\mM_{\mathrm{design}}^\dagger$ denotes the Moore--Penrose
pseudoinverse of $\mM_{\mathrm{design}}$.
\STATE Recover the principal-axis EMT components:
\[
    \widehat{\vv}
    \gets
    \mM_c^{-1}
    \bigl(
      \widehat M^{(1)},
      \widehat M^{(2)},
      \widehat M^{(3)},
      \widehat M^{(4)}
    \bigr)^\top,\quad \widehat{\vv}
    =
    \bigl(
      \widehat m_{1111},
      \widehat m_{2222},
      \widehat m_{1122},
      \widehat m_{1212}
    \bigr)^\top.
\]

\STATE Recover
$(\widehat{\lambda}',\widehat{\mu}',
\widehat{|B|},\widehat m)$
by solving the nonlinear least-squares problem \eqref{eq:nls}.

\STATE Recover the density:
$\widehat{\rho}'
\gets\rho-\widehat M^{(0)}/\widehat{|B|}$.

\STATE Recover two semi-axes: $\widehat a
    =
    \sqrt{
      \frac{(1+\widehat m)\widehat{|B|}}
           {(1-\widehat m)\pi}
    },
    \quad
    \widehat b
    =
    \sqrt{
      \frac{(1-\widehat m)\widehat{|B|}}
           {(1+\widehat m)\pi}}$.

\STATE Resolving the $\pi/2$ ambiguity of $\widehat{\phi}$.

\STATE \textbf{Return}
$(\widehat{\vz},\widehat{\phi},
\widehat{\lambda}',\widehat{\mu}',
\widehat{|B|},\widehat m,\widehat{\rho}')$.

\end{algorithmic}
\end{algorithm}
\subsection{Disk-shaped inclusion}\label{subsec:spherical}

The case of a disk-shaped inclusion differs substantially from that of an elliptic inclusion. Although the EMT can still be uniquely recovered from the PP mode $\fI_{P,P}^{\mathrm{meas}}(\vz^s)$, the three unknown quantities $\lambda'$, $\mu'$, and $|B|$ cannot be uniquely determined from the two independent components of the EMT. Consequently, the corresponding inverse problem is underdetermined.

Nevertheless, useful quantitative information can still be extracted. In particular, we consider the estimation of the inclusion volume. Information on the volume is available from both the density contrast and the elastic moment tensor. Since the density contribution follows directly from the estimate of
$
M^{(0)}=(\rho-\rho')|B|,
$
we focus on the information provided by the EMT. We first consider the case where no prior information on the inclusion is available. In this setting, the following theorem provides constraints on the inclusion volume.

\begin{thm}\label{thm:volume}
    The volume of the inclusion $|B|$ satisfies the following constraints:
    \begin{align}
        &|B|>\max \bra{\frac{m_{1212}\tau}{\mu(\tau+1)},-\frac{m_{1212}}{\mu(\tau+1)}},\\
        & m_{1122}m_{1212}(|B|-B_1)(|B|-B_2)>0,\\
        & (|B|-B_2)\bra{|B|-B_2+\frac{m_{1122}+m_{1212}}{\lambda+\mu}}>0,
    \end{align}
    where
    \begin{align}
        B_1 = -\frac{m_{1212}}{\lambda+2\mu}
        +\frac{m_{1212}\tau}{\mu(\tau+1)}
        +\frac{m_{1212}^2}{m_{1122}\mu(\tau+1)},
        \qquad
        B_2=\frac{m_{1122}}{\lambda+2\mu}
        +\frac{m_{1212}(\tau-1)}{\mu(\tau+1)}.
    \end{align}
\end{thm}

\begin{proof}
See \ref{proof:elliptic_thm}.
\end{proof}

The above constraints are derived without assuming any prior information about the inclusion. In practice, however, it is often reasonable to assume that the Lam\'e parameters and the density of the inclusion lie within prescribed ranges relative to those of the background medium. Such prior information can be incorporated to further improve the estimates of the inclusion volume. Since the contribution associated with the density contrast can be obtained directly from the estimate of $M^{(0)}$, we focus below on improving the bounds derived from the EMT.

\begin{prop}\label{prop:volume}
Assume that
\[
\alpha^-\lambda<\lambda'<\alpha^+\lambda,
\qquad
\beta^-\mu<\mu'<\beta^+\mu,
\]
for some constants $\alpha^\pm,\beta^\pm>0$. Since
\[
(\lambda-\lambda')(\mu-\mu')>0,
\]
there are two admissible cases:
\begin{enumerate}
    \item
    $
    0<\alpha^-<\alpha^+<1,
    \quad
    0<\beta^-<\beta^+<1;
    $
    \item
    $
    1<\alpha^-<\alpha^+,
    \quad
    1<\beta^-<\beta^+.
    $
\end{enumerate}
Then the inclusion volume satisfies
\begin{align}
        &(\alpha^--1)\lambda
        <
        \frac{m_{1122}m_{1212}\mu(\tau+1)(|B|-B_1)}
        {(|B|-B_2)(|B|-B_3)}
        <
        (\alpha^+-1)\lambda,\\
        &|B|\in
        \bra{
        \frac{m_{1212}}{\mu}
        \bra{\frac{\tau}{\tau+1}-\frac{1}{1-\beta^-}},
        \frac{m_{1212}}{\mu}
        \bra{\frac{\tau}{\tau+1}-\frac{1}{1-\beta^+}}
        },
        \qquad
        0<\beta^-<\beta^+<1,\\
        &|B|\in
        \bra{
        \frac{m_{1212}}{\mu}
        \bra{\frac{\tau}{\tau+1}-\frac{1}{1-\beta^+}},
        \frac{m_{1212}}{\mu}
        \bra{\frac{\tau}{\tau+1}-\frac{1}{1-\beta^-}}
        },
        \qquad
        1<\beta^-<\beta^+,
\end{align}
where $B_1$ and $B_2$ are given in Theorem~\ref{thm:volume}, and
\[
B_3=\frac{m_{1212}\tau}{\mu(\tau+1)}.
\]
\end{prop}

\begin{proof}
See \ref{proof_elliptic_prop}.
\end{proof}

For a disk-shaped inclusion, the reconstruction procedure begins by recovering the EMT from the PP mode using the same least-squares procedure described in the previous subsection. The recovered EMT is then used together with the estimate of $M^{(0)}$ to infer the inclusion volume. The EMT provides the bounds in Theorem~\ref{thm:volume}, while the estimate of $M^{(0)}$ supplies the information associated with the density contrast. When prior information on the material parameters is available, the sharper bounds in Proposition~\ref{prop:volume} can instead be employed.

We conclude this section with a remark on inclusions of more general
shape. Although the identification procedures developed above are
formulated for elliptic and disk-shaped inclusions, the leading-order
scattered field in the small-volume regime depends on the inclusion
geometry only through its EMT, together with the density--volume
contrast. A general small inclusion can therefore be represented, at
leading-order, by an equivalent ellipse having the same EMT. The
proposed procedure then recovers the material and geometric parameters
of this equivalent ellipse rather than the exact shape of the original
inclusion. Thus, the method remains applicable to small inclusions of
more general shape at the level of their leading-order effective
representation.

\section{Numerical experiments}\label{sec:numerics}
We now present numerical experiments illustrating the imaging and reconstruction procedures. All simulations are performed in MATLAB on an Apple M2 processor with an 8-core CPU. The numerical setup is as follows.

We assume that the inclusion is made of iron, while the background medium is fused silica. Accordingly, the material parameters are given by
$\rho = 2.20~\mathrm{g/cm}^3$, $c_P^{(\text{silica})} = 5.97\times10^5~\mathrm{cm/s}$, $c_S^{(\text{silica})} = 3.76\times10^5~\mathrm{cm/s}$, $\rho' = 7.54~\mathrm{g/cm}^3$, $c_P^{(\text{iron})} = 5.80\times10^5~\mathrm{cm/s}$, and $c_S^{(\text{iron})} = 3.10\times10^5~\mathrm{cm/s}$. The corresponding Lamé parameters are
\[
\lambda = 1.62\times10^{11}~\mathrm{dyn/cm}^2,\quad
\mu = 3.11\times10^{11}~\mathrm{dyn/cm}^2,\quad
\lambda' = 1.09\times10^{12}~\mathrm{dyn/cm}^2,\quad
\mu' = 7.25\times10^{11}~\mathrm{dyn/cm}^2.
\]

The operating frequency is set to $2$ MHz. The corresponding typical $P$ and $S$ wavelengths in silica are approximately \(0.2985\,\mathrm{cm}\) and \(0.188\,\mathrm{cm}\), respectively. The inclusion is modeled as an ellipse with semi-axes $a = 0.006~\mathrm{cm}$ and $b = 0.0015~\mathrm{cm}$ . A total of 8192 collocated sources and receivers are uniformly distributed on a circle of radius $R = 4000~\mathrm{cm}$. The inclusion is assumed to be located at $(400~\mathrm{cm}, 400~\mathrm{cm})$.

In Section \ref{subsec:ex_noiseless}, the scattered field is generated using
\begin{align}\label{eq:exp_data_gen}
\vu^{(j)}_{\mathrm{s},\mathrm{meas}}(\vx_l,\vx_k)
= \vu^{(j)}_{\mathrm{scat}}(\vx_l,\vx_k)
+ \vw^{(j)}(\vx_l,\vx_k),
\end{align}
where $\vw^{(j)}$ denotes the additive noise term. Here, $\vu^{(j)}_{\mathrm{scat}}(\vx_l,\vx_k)$ is computed using the leading-order small-volume expansion. The validity of this expansion in the regime we consider follows from the standard
small-volume expansion theory for elastic inclusions; see
\cite{3d546238-6f8a-3f6d-b0b4-3883234b540d,
ammari2004reconstruction,doi:10.1137/070696076}.
For a fixed operating frequency and source and receiver points uniformly
separated from the inclusion, the leading-order contribution is of
order $\varepsilon^2$ in two dimensions, while the remainder is of
order $\mathcal O(\varepsilon^3)$. For the physical parameters used below,
\[
    \frac{a}{c_P/f}\approx 2.01\times10^{-2},
    \qquad
    \frac{a}{c_S/f}\approx 3.19\times10^{-2},
\]
so the inclusion is small relative to both elastic wavelengths. These
ratios, together with the cited remainder estimate, support the use of
the leading-order small-volume approximation.

In Section \ref{subsec:full-wave}, the scattered field $\vu_{\mathrm{s},\mathrm{meas}}(\vx_l,\vx_k)=\vu_\varepsilon-\vu_{\mathrm{inc}}$ is computed by solving the elastic transmission problem \eqref{eq:scattering_problem} using a boundary collocation method that combines the method of fundamental solutions (MFS) \cite{FAIRWEATHER2003759} with regular cylindrical wave expansions. The exterior scattered field is represented by outgoing elastic fundamental solutions with fictitious sources inside the inclusion, while the interior total field is expanded in regular compressional and shear cylindrical waves. These representations satisfy the respective Navier equations, and the exterior basis incorporates the radiation condition. The expansion coefficients are determined by enforcing continuity of displacement and traction on the interface. We use 384 fictitious sources and 1536 boundary collocation points. The incident field expansion is truncated at order 12, while the interior expansion retains angular modes up to order 20. After column scaling, the resulting least-squares system is solved by truncated singular value decomposition with a relative truncation threshold of $10^{-13}$.

In both cases, the scattered field is sampled at the receiver
locations $\{\vx_l\}_{l=1}^{N_r}$ to obtain the measured data
$\vu_{\mathrm{s},\mathrm{meas}}(\vx_l,\vx_k)$.
Throughout this section, we refer to the data generated using the
small-volume expansion as \emph{asymptotic data} and those obtained
by solving the transmission problem as \emph{full-wave data}.

We consider the imaging functional defined by
\begin{align}
\fI_{\alpha,\beta}^{\mathrm{meas}}(\vz^s)
= \sum_{j=1}^2 \int_{\partial\Omega} e^{-i\kappa_\beta|\vy-\vz^s|}\ve_j \cdot
\int_{\partial\Omega} e^{-i\kappa_\alpha|\vx-\vz^s|}\vu^{(j)}_{\mathrm{s},\mathrm{meas}}(\vx,\vy)\, d\sigma(\vx)\, d\sigma(\vy),
\quad \alpha,\beta = P,S.
\end{align}

The relative error of an estimate $\hat{\theta}$ with respect to the ground truth $\theta^*$ is defined as
\[
E_{\mathrm{relative}} = \frac{|\hat{\theta}-\theta^*|}{|\theta^*|}.
\]
\subsection{Imaging and identification with noiseless asymptotic data}\label{subsec:ex_noiseless}

In this section, we present numerical results for inclusion imaging and
the identification of its unknown material and geometric parameters using the asymptotic data.

\subsubsection{Imaging the inclusion}

The synthetic data are generated according to
\eqref{eq:exp_data_gen}, with the additive noise set to zero,
$\vw^{(j)}=0$.
For imaging, we consider the search domain
\[
[400-1/3,400+1/3]^2,
\]
uniformly discretized into a $101\times101$ grid.

\begin{figure}
    \centering
    \includegraphics[width=\linewidth]{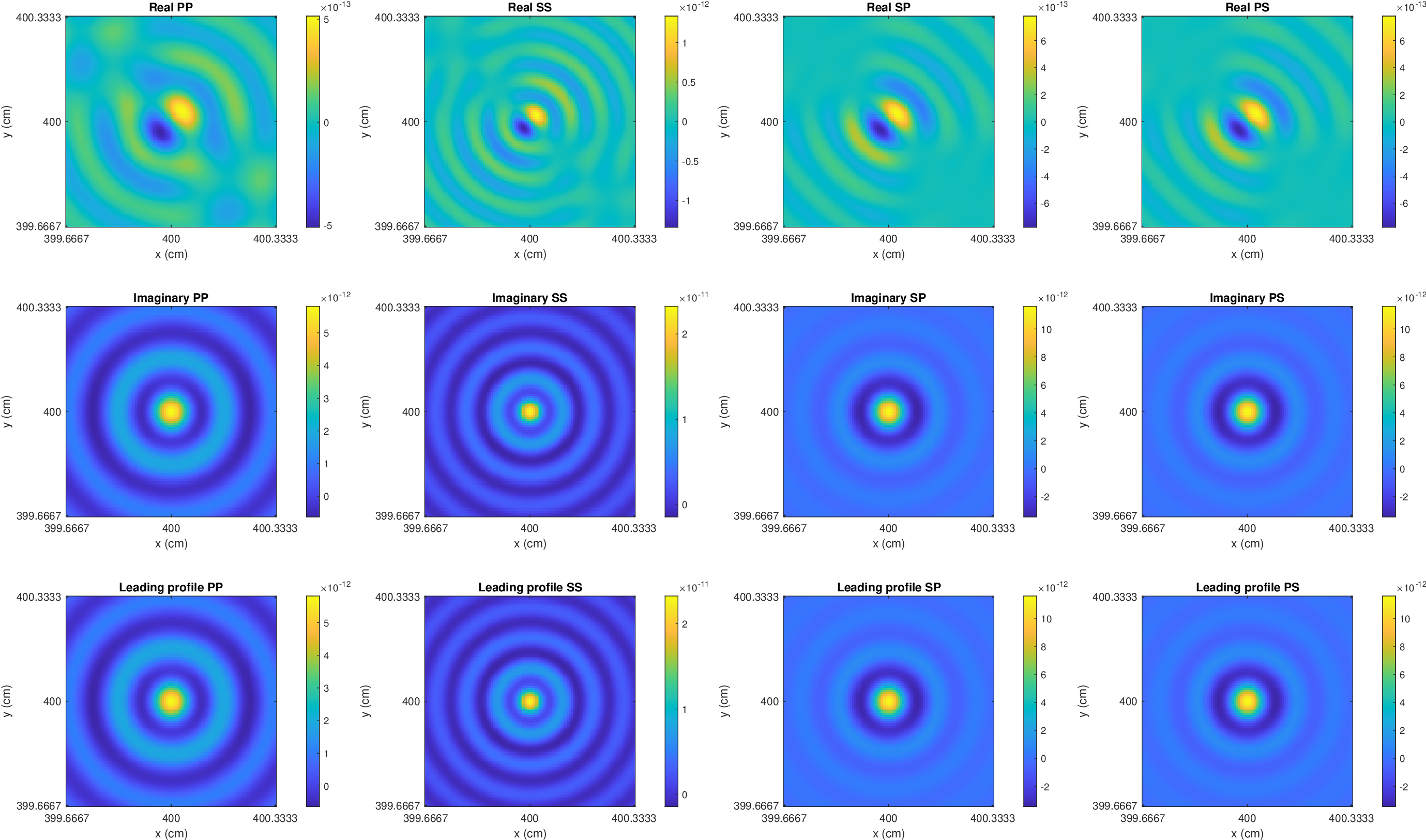}
    \caption{Comparison of the numerically evaluated imaging functionals
and their leading-order asymptotic predictions for the four
polarization modes. The excitation frequency is
\(f=2\,\mathrm{MHz}\), and the inclusion is centered at
\(\vz=(400\,\mathrm{cm},400\,\mathrm{cm})\) with orientation $\phi=0$. The first and second rows show,
respectively, the real and imaginary parts of the numerical imaging
functionals, while the third row shows the imaginary parts of the
corresponding leading-order asymptotic profiles
\(\fI^{\mathrm{asym}}_{\alpha,\beta}\) defined in
\eqref{eq:I_asym}; their real parts vanish. From left to right, the
columns correspond to the PP, SS, SP, and PS modes.}
    \label{fig:compare}
\end{figure}

Figure~\ref{fig:compare} compares the numerically evaluated imaging
functionals with their theoretical asymptotic profiles for the four
polarization modes. Both the real and imaginary parts of the numerical
results are displayed, while localization is performed using the
imaginary part, in accordance with the theoretical prediction that the
leading-order imaging functionals are purely imaginary.

\subsubsection{Discussion of the model error}
\label{subsec:model_error}
We first recall the origin of the cross-mode contributions that are
neglected in deriving the explicit leading-order imaging profiles. The
elastic Green tensor contains both compressional ($P$) and shear ($S$)
wave components. When the measured field is back-propagated using a
prescribed polarization mode, the Green-tensor component with the same
polarization is phase matched, whereas the component with the other
polarization generates a cross-mode contribution. Under the
scale-separation condition
\eqref{eq:large-frequency-distance}, the cross-mode contributions are
asymptotically negligible and are therefore omitted from the explicit
leading-order profiles. When this condition is not sufficiently
satisfied, however, these neglected terms may become significant.

In this subsection, the term \emph{model error} refers specifically to
the discrepancy resulting from the omission of these cross-mode
contributions. To distinguish this error from measurement noise, the
truncation error of the small-volume expansion, and the far-field
approximation error of the Green tensor, we introduce a reference
imaging functional
$\fI_{\alpha,\beta}^{\mathrm{ref}}$. This functional is evaluated from
the noiseless leading-order small-volume model using the same
far-field approximation of the Green tensor as in the asymptotic
derivation, but with all cross-mode contributions retained. 

For each polarization mode $(\alpha,\beta)$, we define the pointwise
scale-separation model error by
\begin{equation}\label{eq:model_error_pointwise}
    \mathcal E_{\alpha,\beta}^{\mathrm{model}}(\vz^s)
    :=
    \fI_{\alpha,\beta}^{\mathrm{ref}}(\vz^s)
    -
    \fI_{\alpha,\beta}^{\mathrm{asym}}(\vz^s),
\end{equation}
where $\fI_{\alpha,\beta}^{\mathrm{asym}}$ denotes the explicit
leading-order profile defined in \eqref{eq:I_asym}. Thus,
$\mathcal E_{\alpha,\beta}^{\mathrm{model}}$ isolates the effect of the
cross-mode contributions that are negligible under
\eqref{eq:large-frequency-distance} but may become significant outside
this asymptotic regime.

For noisy measured or synthetic data, the total discrepancy from the
asymptotic profile can be decomposed as
\begin{equation}\label{eq:total_error_decomposition}
\begin{aligned}
    \fI_{\alpha,\beta}^{\mathrm{meas}}(\vz^s)
    -
    \fI_{\alpha,\beta}^{\mathrm{asym}}(\vz^s)
    &=
    \Bigl[
        \fI_{\alpha,\beta}^{\mathrm{meas}}(\vz^s)
        -
        \fI_{\alpha,\beta}^{\mathrm{ref}}(\vz^s)
    \Bigr]
    +
    \mathcal E_{\alpha,\beta}^{\mathrm{model}}(\vz^s).
\end{aligned}
\end{equation}
The bracketed term accounts for measurement noise and any additional
discrepancies between the data and the reference model, whereas
$\mathcal E_{\alpha,\beta}^{\mathrm{model}}$ represents the
scale-separation model error considered in this subsection.

To isolate this effect, we keep the excitation frequency, inclusion
size, and measurement radius fixed and move the inclusion closer to
the origin. The excitation frequency is
\[
    f=2\,\mathrm{MHz},
    \qquad
    \omega=2\pi f,
\]
and the elliptic inclusion has the same semi-axes as in the preceding
example:
\[
    a=0.006\,\mathrm{cm},
    \qquad
    b=0.0015\,\mathrm{cm}.
\]
Its center is placed at
\[
    \vz=(0.25\,\mathrm{cm},0.25\,\mathrm{cm}).
\]
For these parameters,
\[
    \frac{\omega|\vz|}{2\pi c_P}\approx1.18,
    \qquad
    \frac{\omega|\vz|}{2\pi c_S}\approx1.88.
\]
Thus, the distance of the inclusion from the origin is only about
$1.18$ $P$ wavelengths and $1.88$ $S$ wavelengths, so
condition~\eqref{eq:large-frequency-distance} is not well satisfied.

By contrast, conditions~\eqref{eq:small-volume-regime} and
\eqref{eq:geometric-separation} remain valid. In particular, the
inclusion remains small relative to both elastic wavelengths, while
the source-to-inclusion and inclusion-to-receiver distances remain
approximately equal to $R$. Hence, the arguments of the Hankel
functions associated with wave propagation remain large, and the
corresponding large-argument approximations of the Green tensor remain
accurate. The principal change in this experiment is therefore the
violation of condition~\eqref{eq:large-frequency-distance}, rather
than a change in the small-volume or Green-tensor far-field regimes.

\begin{figure}
    \centering
    \includegraphics[width=\linewidth]{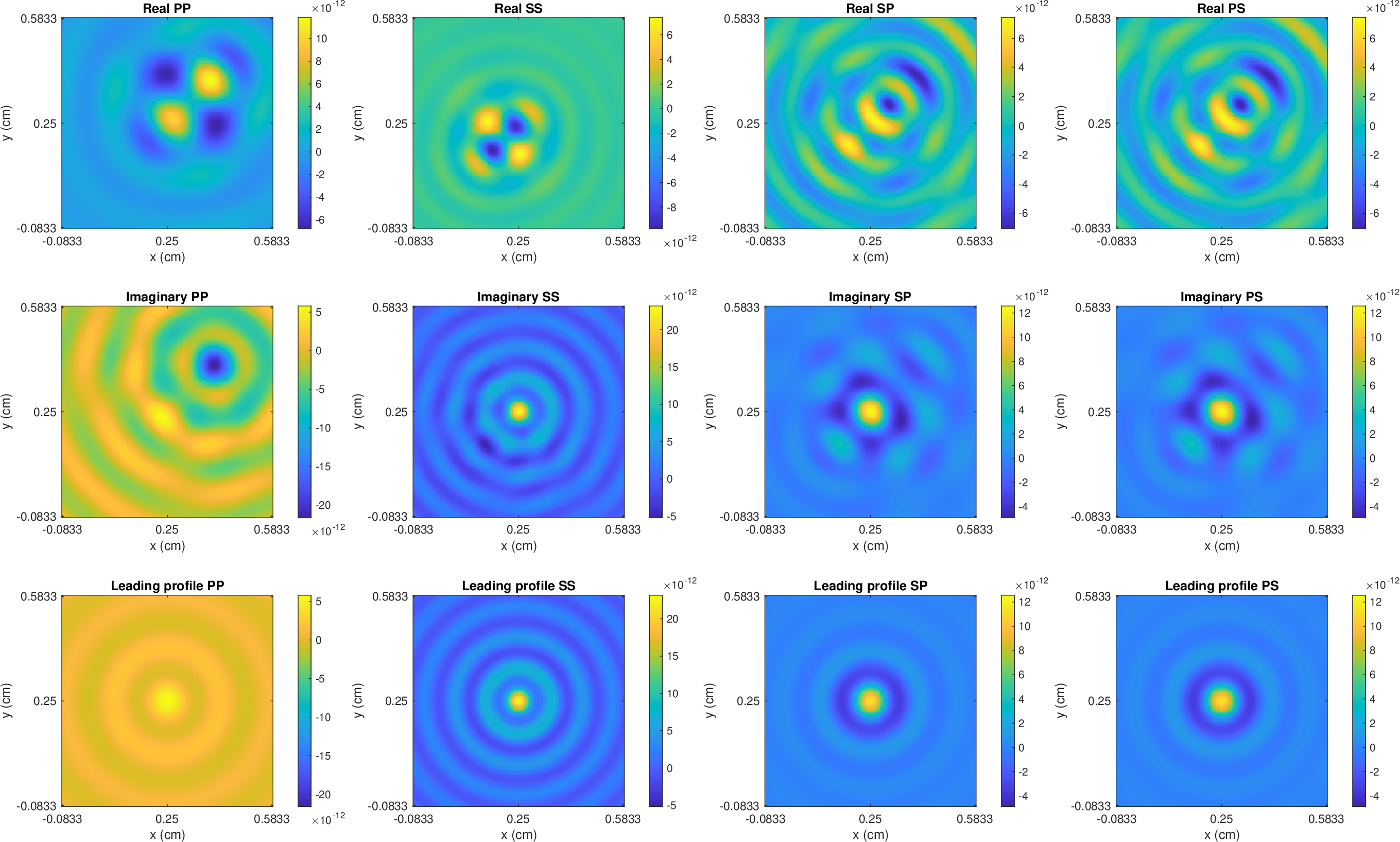}
    \caption{Comparison of the numerically evaluated imaging functionals
    and their leading-order asymptotic predictions for the four
    polarization modes. The excitation frequency is
    \(f=2\,\mathrm{MHz}\), and the inclusion is centered at
    \(\vz=(0.25\,\mathrm{cm},0.25\,\mathrm{cm})\) with orientation $\phi=0$. The first and
    second rows show, respectively, the real and imaginary parts of the
    numerical imaging functionals, while the third row shows the
    imaginary parts of the corresponding leading-order asymptotic
    profiles \(\fI_{\alpha,\beta}^{\mathrm{asym}}\) defined in
    \eqref{eq:I_asym}; their real parts vanish. From left to right, the
    columns correspond to the PP, SS, SP, and PS modes.}
    \label{fig:compare_6e6}
\end{figure}

Figure~\ref{fig:compare_6e6} compares the numerically evaluated imaging
functionals with their corresponding leading-order asymptotic
predictions.

Compared with Figure~\ref{fig:compare}, the discrepancy between the
numerically evaluated imaging functionals and their theoretical
predictions is considerably more pronounced. In particular, the
numerical imaging functionals exhibit non-negligible real parts, which
are absent from the leading-order asymptotic formulas. The imaging
quality also deteriorates: the maxima of the magnitudes of their imaginary parts may
be shifted away from the true inclusion location, and their spatial
profiles no longer closely match the theoretical predictions.

Since the frequency, inclusion size, and measurement radius are kept
fixed, the small-volume parameters and the large arguments of the
Hankel functions remain essentially unchanged. The increased
discrepancy is therefore attributed primarily to the cross-mode
contributions neglected under
condition~\eqref{eq:large-frequency-distance}.

A complete characterization of the model error is beyond the scope of
the present paper. Instead, we illustrate one of its principal
mechanisms by considering the density-contrast contribution to the
imaging functional. This contribution already captures the effect of
the cross-mode terms neglected under
condition~\eqref{eq:large-frequency-distance} and allows their
dependence on the wavenumbers, inclusion location, and Lam\'e
parameters to be made explicit.

The elasticity-contrast contribution can be analyzed in a similar
manner. Its derivation involves additional derivatives of the Green
tensor and contractions with the EMT and is therefore algebraically
more involved. Nevertheless, it contains analogous phase-matched and
cross-mode terms governed by the same oscillatory phase mechanism. We
therefore restrict the following calculation to the density-contrast
contribution; see Section~\ref{sec:psf_resolution} for the corresponding
decomposition of the imaging functional.

For the $P$ wave backpropagation factor, consider
\begin{align}
    \mF_P(\vz^s,\vz)
    &:=
    \int_{\partial\Omega}
    \mGamma^\omega(\vx-\vz)
    e^{-i\kappa_P|\vz^s-\vx|}
    \,d\sigma(\vx)
    \notag\\
    &\sim
    \int_{\partial\Omega}
    \sqrt{\frac{2}{\pi\kappa_S r_{xz}}}
    \frac{\mI-\widehat{\vr}_{xz}\widehat{\vr}_{xz}^{\top}}
         {4\mu}
    e^{i\left(
        \kappa_S|\vx-\vz|
        -
        \kappa_P|\vz^s-\vx|
    \right)}
    \,d\sigma(\vx)
    \notag\\
    &\quad+
    \int_{\partial\Omega}
    \sqrt{\frac{2}{\pi\kappa_P r_{xz}}}
    \frac{\widehat{\vr}_{xz}\widehat{\vr}_{xz}^{\top}}
         {4(\lambda+2\mu)}
    e^{i\kappa_P\left(
        |\vx-\vz|-|\vz^s-\vx|
    \right)}
    \,d\sigma(\vx)
    \notag\\
    &\triangleq
    \mF_{P,\mathrm{cross}}
    +
    \mF_{P,\mathrm{match}},
    \label{eq:FP_decomposition}
\end{align}
where
\[
    \vr_{xz}:=\vx-\vz,
    \qquad
    r_{xz}:=|\vr_{xz}|,
    \qquad
    \widehat{\vr}_{xz}:=\frac{\vr_{xz}}{r_{xz}}.
\]
Here, $\mF_{P,\mathrm{match}}$ is the phase-matched $P$ wave
contribution, whereas $\mF_{P,\mathrm{cross}}$ is the cross-mode
contribution arising from the $S$ wave component of the Green tensor.

To estimate the latter, we consider sampling points in a neighborhood
of the inclusion and write
\[
    \vz^s=\vz+\varepsilon^{1/2}\vr^s.
\]
For fixed frequency, the cross-mode phase satisfies
\begin{align}
    \kappa_S|\vx-\vz|
    -
    \kappa_P|\vz^s-\vx|
    &=
    (\kappa_S-\kappa_P)|\vx-\vz|
    +
    \mathcal{O}(\varepsilon^{1/2})
    \notag\\
    &\approx
    \Delta\kappa
    \left(
        R-\langle\widehat{\vx},\vz\rangle
    \right),
    \label{eq:cross_phase_P}
\end{align}
where
\[
    \Delta\kappa:=\kappa_S-\kappa_P>0,
    \qquad
    \widehat{\vx}:=\frac{\vx}{|\vx|}.
\]
Applying the stationary phase method to
\eqref{eq:FP_decomposition} gives
\[
    \mF_{P,\mathrm{cross}}
    =
    \mathcal{O}\left(
        \frac{\sqrt{R}}
        {\mu\sqrt{\kappa_S}\sqrt{\Delta\kappa|\vz|}}
    \right).
\]
On the other hand, the phase-matched contribution satisfies
\[
    \mF_{P,\mathrm{match}}
    =
    \mathcal{O}\left(
        \frac{1}{\lambda+2\mu}
        \sqrt{\frac{R}{\kappa_P}}
    \right),
\]
as shown in \eqref{eq:F_p1}. Consequently,
\begin{equation}\label{eq:relative_cross_P}
    \left|
        \frac{\mF_{P,\mathrm{cross}}}
             {\mF_{P,\mathrm{match}}}
    \right|
    =
    \mathcal{O}\left(
        \left(
            \frac{\lambda+2\mu}{\mu}
        \right)^{3/4}
        \frac{1}{\sqrt{\Delta\kappa|\vz|}}
    \right).
\end{equation}

A similar decomposition holds for the $S$ wave backpropagation factor:
\begin{align}
    \mF_S(\vz^s,\vz)
    &:=
    \int_{\partial\Omega}
    \mGamma^\omega(\vx-\vz)
    e^{-i\kappa_S|\vz^s-\vx|}
    \,d\sigma(\vx)
    \notag\\
    &\sim
    \int_{\partial\Omega}
    \sqrt{\frac{2}{\pi\kappa_P r_{xz}}}
    \frac{\widehat{\vr}_{xz}\widehat{\vr}_{xz}^{\top}}
         {4(\lambda+2\mu)}
    e^{i\left(
        \kappa_P|\vx-\vz|
        -
        \kappa_S|\vz^s-\vx|
    \right)}
    \,d\sigma(\vx)
    \notag\\
    &\quad+
    \int_{\partial\Omega}
    \sqrt{\frac{2}{\pi\kappa_S r_{xz}}}
    \frac{\mI-\widehat{\vr}_{xz}\widehat{\vr}_{xz}^{\top}}
         {4\mu}
    e^{i\kappa_S\left(
        |\vx-\vz|-|\vz^s-\vx|
    \right)}
    \,d\sigma(\vx)
    \notag\\
    &\triangleq
    \mF_{S,\mathrm{cross}}
    +
    \mF_{S,\mathrm{match}}.
    \label{eq:FS_decomposition}
\end{align}
Here, the first term is the cross-mode $P$ wave contribution, whereas
the second is the phase-matched $S$ wave contribution. Repeating the
stationary phase argument yields
\[
    \mF_{S,\mathrm{cross}}
    =
    \mathcal{O}\left(
        \frac{\sqrt{R}}
        {(\lambda+2\mu)\sqrt{\kappa_P}
         \sqrt{\Delta\kappa|\vz|}}
    \right)
\]
and
\[
    \mF_{S,\mathrm{match}}
    =
    \mathcal{O}\left(
        \frac{1}{\mu}
        \sqrt{\frac{R}{\kappa_S}}
    \right).
\]
It follows that
\begin{equation}\label{eq:relative_cross_S}
    \left|
        \frac{\mF_{S,\mathrm{cross}}}
             {\mF_{S,\mathrm{match}}}
    \right|
    =
    \mathcal{O}\left(
        \left(
            \frac{\mu}{\lambda+2\mu}
        \right)^{3/4}
        \frac{1}{\sqrt{\Delta\kappa|\vz|}}
    \right).
\end{equation}

The estimates above are derived for the density-contrast contribution.
For the elasticity contrast, the derivatives of the Green tensor give
rise to additional amplitude factors and tensor contractions, but the
cross-mode terms contain the same phase mismatch between the
$P$ and $S$ wave components. A stationary phase analysis therefore
leads to the same qualitative decay with respect to
$\Delta\kappa|\vz|$, although the corresponding prefactors depend more
intricately on the EMT and the Lam\'e parameters.

Consequently, cross-mode residual terms of this type occur in all four
polarization-dependent imaging functionals for both density and
elasticity contrasts. Equations~\eqref{eq:relative_cross_P} and
\eqref{eq:relative_cross_S} make this mechanism explicit for the
density-contrast contribution and show that its magnitude decreases as
the dimensionless parameter $\Delta\kappa|\vz|$ increases.

Moreover, since $\lambda+2\mu>\mu$, the Lam\'e-parameter factors in
\eqref{eq:relative_cross_P} and \eqref{eq:relative_cross_S} show that
the relative cross-mode contribution is larger for $P$ wave
backpropagation than for $S$ wave backpropagation. For the density-contrast contribution, this explains why the PP mode
is most strongly affected, the SS mode is least affected, and the
mixed PS and SP modes exhibit intermediate behavior. The
elasticity-contrast contribution exhibits the same qualitative trend,
although its precise magnitude also depends on the EMT components.

\subsubsection{Identification of the properties of the inclusion with $\phi=\pi/6$}\label{subsubsec:exp_phi_30}
\begin{figure}
    \centering
    \includegraphics[width=0.9\linewidth]{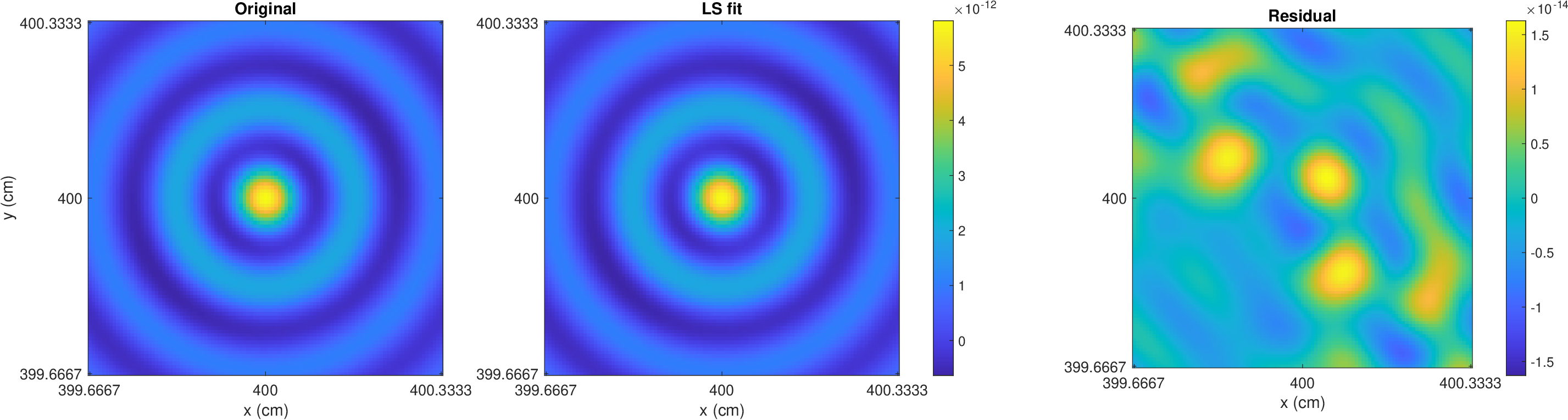}
    \caption{Results of the least-squares fitting experiment. From left to right, the panels show the imaginary part of the original image, the image reconstructed by least-squares fitting, and the corresponding residual.}
    \label{fig:LS_experiments}
\end{figure}

After utilizing the $SS$ mode imaging functional to reconstruct the location of the inclusion, we follow the algorithm to estimate the orientation. The numerical result shows that the estimated orientation is $\hat{\phi} = 0.5221 ~\mathrm{rad}$, with $0.2914\%$ relative error.
We then recover the coefficients $\{\widehat{M}^{(j)}\}_{j=0}^{4}$ by least-squares fitting. The image is evaluated on a uniform grid with spacing $1/150~\mathrm{cm}$ over the square $[-1/3,1/3]^2$, translated to be centered at the inclusion. The relative fitting error over the entire imaging grid is
\[
    \frac{\|I_{\mathrm{orig}}-I_{\mathrm{LS}}\|_F}
         {\|I_{\mathrm{orig}}\|_F}
    = 0.48\%.
\]
Here, $I_{\mathrm{orig}}$ denotes the imaginary part of the
numerically evaluated PP imaging functional, and $I_{\mathrm{LS}}$
denotes its least-squares approximation.
Both are represented as two-dimensional matrices, and $\|\cdot\|_F$ denotes the Frobenius norm. The corresponding results—including the original image, least-squares reconstruction, and residual—are shown in Figure~\ref{fig:LS_experiments}.

To solve the resulting nonlinear least-squares problem, we adopt a multi-start strategy to improve robustness. Specifically, 50 random initializations are used. In each run, the initial guesses for $\lambda'$, $\mu'$, $|B|$, and $m$ are perturbed within $\pm 40\%$ of their true values. The optimization problem is solved from each initial guess using MATLAB’s \texttt{fminunc} function with a quasi-Newton method, where a BFGS update is employed to approximate the inverse Hessian. The optimization is unconstrained, with strict convergence tolerances of $10^{-12}$ for both optimality and step-size criteria. The maximum number of iterations and function evaluations is set to 2000 and 20000, respectively. Each starting point is optimized independently, and the estimate with the smallest residual among the 50 starts is retained. The boxplots summarize the relative errors of the 100 retained estimates.

\begin{figure}
    \centering
    \includegraphics[width=0.5\linewidth]{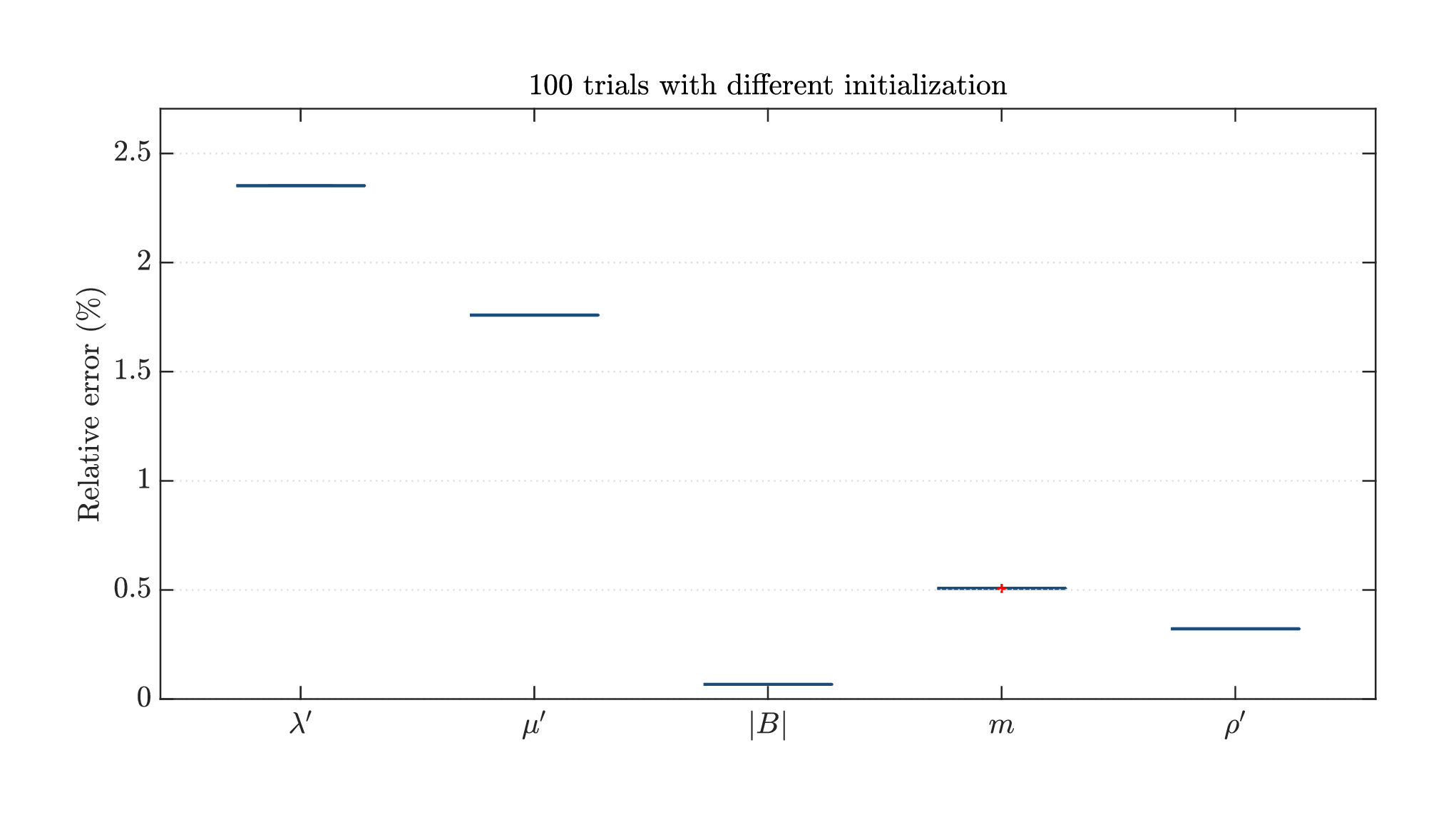}
    \caption{Distribution of the relative errors for different reconstructed parameters over 100 trials using a 50-start optimization strategy.}
    \label{fig:relative_error_mont}
\end{figure}

We observe that the results across different trials are stable under the 50-start strategy. To quantify the variability, we perform 100 Monte Carlo experiments and report the relative errors of the reconstructed parameters in Figure~\ref{fig:relative_error_mont}.

\subsection{Imaging and identification with full-wave data}\label{subsec:full-wave}
In this section, we use the same setup for the inclusion and background medium as in
Section~\ref{subsubsec:exp_phi_30}, and we use the full-wave data for the scattered field.

\subsubsection{Noiseless full-wave data}
In this section, we present and discuss the result for the inclusion identification with noiseless full-wave data. For ease of comparison,
Table~\ref{tab:fullwave-reconstruction} summarizes the reconstruction
results obtained from asymptotic and full-wave data. The larger
reconstruction errors observed for full-wave data are consistent
with the truncation error introduced by the small-volume
approximation. In addition, all reported
relative errors remain below $8.1\%$, indicating that the proposed
method remains effective in this setting despite the use of a
leading-order asymptotic reconstruction model.
\begin{table}[htbp]
\centering
\caption{Comparison of parameter reconstruction from noiseless
asymptotic and full-wave data for an elliptic inclusion with
orientation $\phi=\pi/6$.}
\label{tab:fullwave-reconstruction}

\setlength{\tabcolsep}{4pt}
\begin{tabular}{lccccc}
\toprule
Parameter & True value & Asymptotic & Error (\%) & Full-wave & Error (\%) \\
\midrule
$\phi$ (rad)
& $0.5236$ & $0.5221$ & $0.2914$ & $0.5220$ & $0.3049$ \\

$\lambda^\prime$ (dyn/cm$^2$)
& $1.087\times10^{12}$ & $1.113\times10^{12}$ & $2.352$
& $9.992\times10^{11}$ & $8.097$ \\

$\mu^\prime$ (dyn/cm$^2$)
& $7.246\times10^{11}$ & $7.373\times10^{11}$ & $1.759$
& $6.992\times10^{11}$ & $3.499$ \\

$|B|$ (cm$^2$)
& $2.827\times10^{-5}$ & $2.826\times10^{-5}$ & $0.06743$
& $2.935\times10^{-5}$ & $3.820$ \\

$m$
& $0.6000$ & $0.5970$ & $0.5071$ & $0.6165$ & $2.752$ \\

$\rho^\prime$ (g/cm$^3$)
& $7.540$ & $7.564$ & $0.3217$ & $7.452$ & $1.163$ \\
\bottomrule
\end{tabular}
\end{table}

\subsubsection{Imaging and identification with noisy full-wave data}\label{subsec:ex_noise}
\begin{figure}[htbp]
    \centering
    \includegraphics[width=0.6\linewidth]{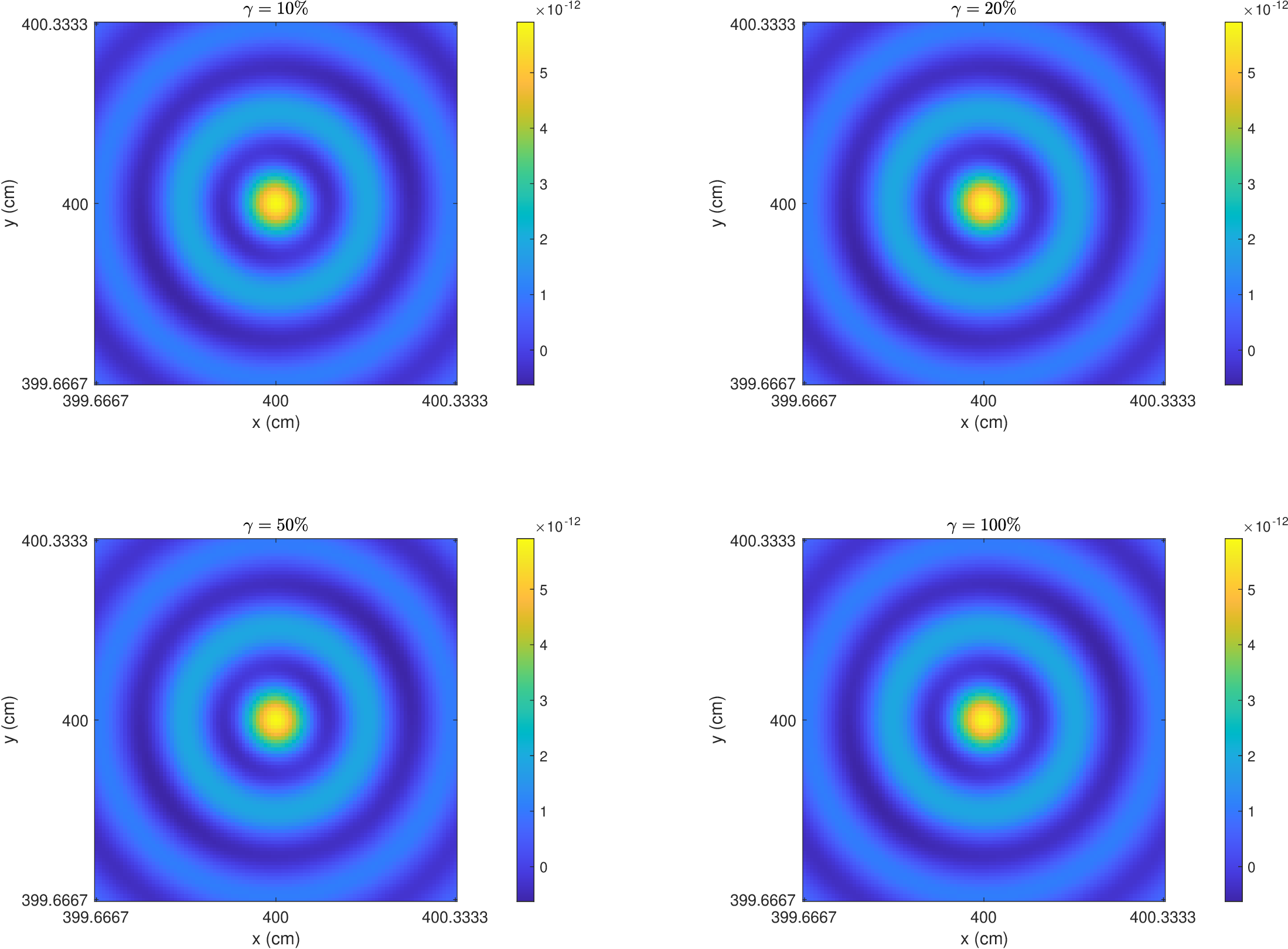}
    \caption{Visualization of PP mode results under different noise levels. The noise level is indicated in each panel. All four images are centered at $(400~ \mathrm{cm},400 ~\mathrm{cm})$.}
    \label{fig:noisy_image} 
\end{figure}

\begin{figure}[htbp]
    \centering
    \includegraphics[width=0.7\linewidth]{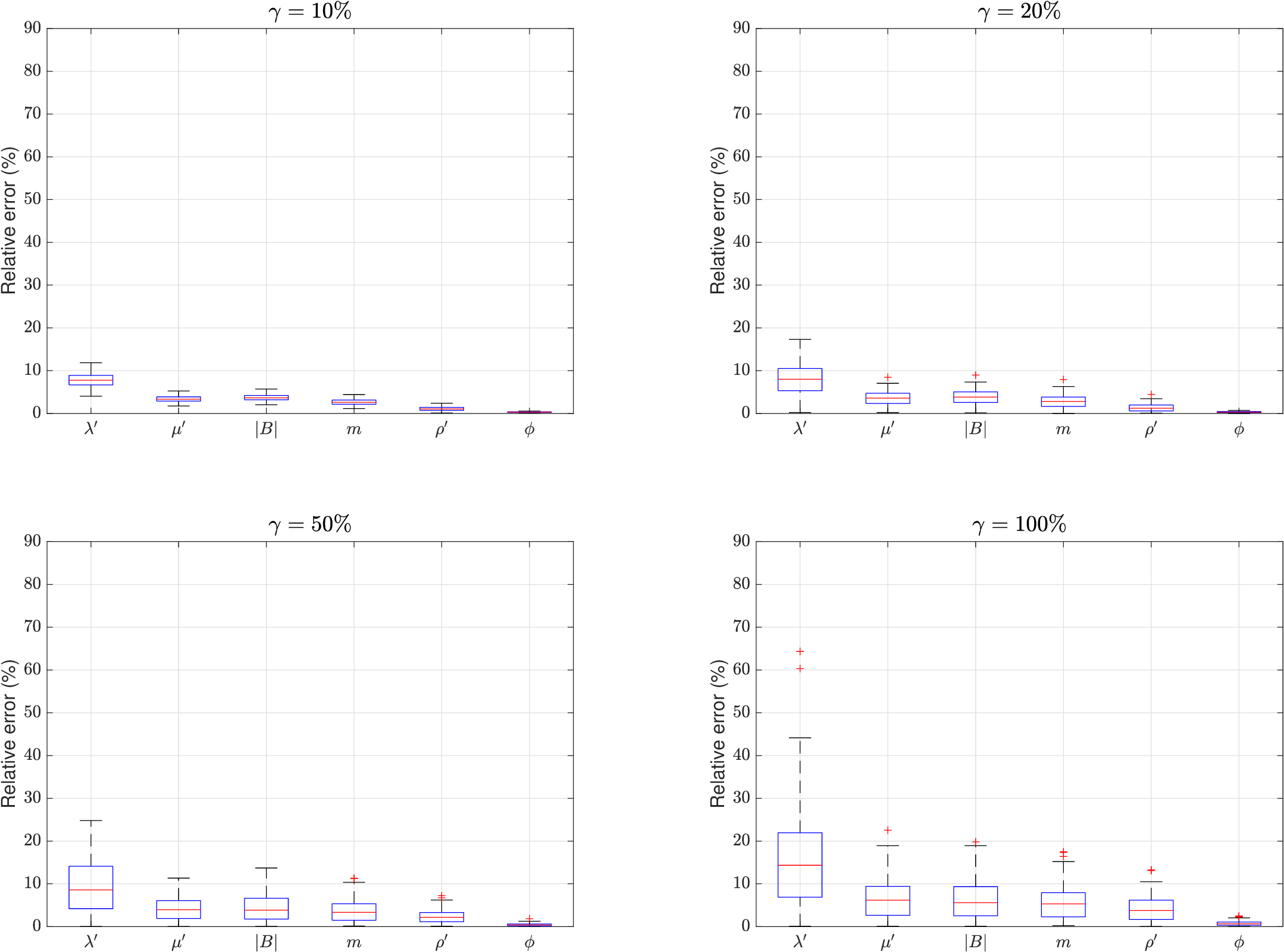}
    \caption{Relative errors of the reconstructed parameters, with the
    $y$-axis fixed to $[0,90\%]$ to facilitate comparison across
    different noise levels. The number of collocated sources and
    receivers is $8192$, and the noise level $\gamma$ is indicated in
    each panel.}
    \label{fig:noisy_relative_error}
\end{figure}

In this section, we consider the reconstruction from noisy data. The measurement noise is modeled as
\[
\vw^{(j)}
=
\gamma\,
\frac{\|\vu_{\mathrm{scat}}^{(j)}\|_F}{\sqrt{N_{\mathrm{num}}}}
\,\vn^{(j)},
\]
where \(N_{\mathrm{num}}\) denotes the total number of entries in the scattered field matrix \(\vu_{\mathrm{scat}}^{(j)}\), \(\|\cdot\|_F\) denotes the Frobenius norm, and \(\gamma>0\) specifies the relative noise level. The random matrix \(\vn^{(j)}\) consists of independent and identically distributed complex Gaussian entries. Specifically, each entry satisfies
\[
n_{l,k}^{(j)} \sim \mathcal{CN}(0,1),
\]
where \(\mathcal{CN}(0,1)\) denotes the standard complex Gaussian distribution.
Figure \ref{fig:noisy_image} shows the imaginary part of PP mode iamging functional for different noise levels. We observe that the imaging of the inclusion using the PP mode is stable under additional noise.
Next, we present Monte Carlo results for the identification of inclusion parameters under different noise levels, as shown in Figure~\ref{fig:noisy_relative_error}. A total of 100 independent experiments are conducted for each noise setting.

We observe that the stability of the parameter identification deteriorates as the level of additive noise increases. For clarity, the $y$-axis is truncated to the range $[0,90\%]$. 
Overall, the proposed method provides reliable parameter identification under moderate noise levels. Even under the most challenging noise conditions, the median relative error for all parameters remains below 17\%, indicating that the method retains useful estimation capability despite severe measurement corruption.

\begin{figure}[htbp]
    \centering
    \includegraphics[width=0.7\linewidth]{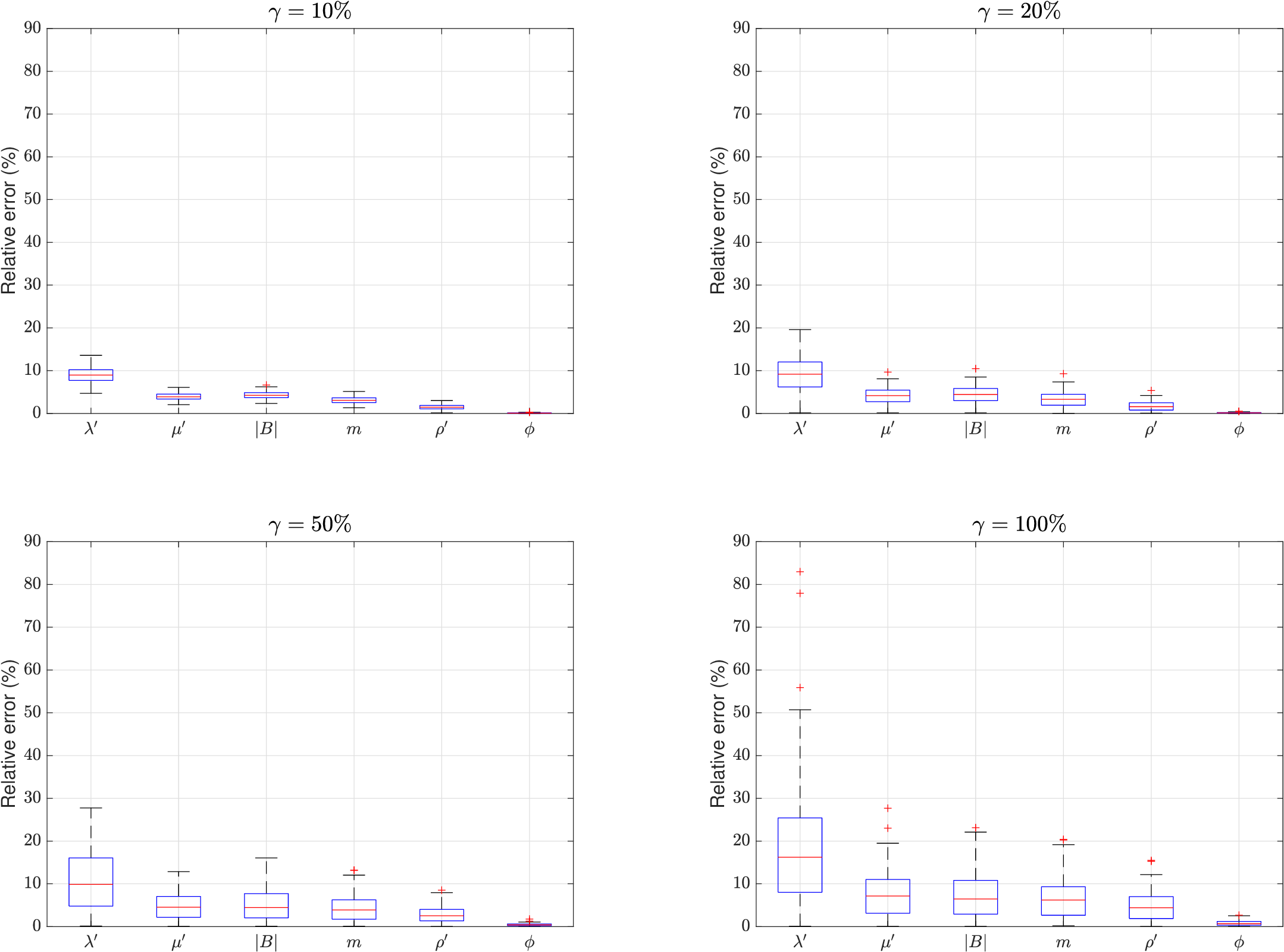}
    \caption{Relative errors of the reconstructed parameters, with the
    $y$-axis fixed to $[0,90\%]$ to facilitate comparison across
    different noise levels. The number of collocated sources and
    receivers is $7050$, and the noise level $\gamma$ is indicated in
    each panel.}
    \label{fig:noisy_relative_error_1k}
\end{figure}

In the following experiment, we reduce the number of collocated sources and receivers to $7050$, while keeping all other experimental settings identical to those in the previous example. The relative errors of the reconstructed parameters are presented in Figure~\ref{fig:noisy_relative_error_1k}. As before, the y-axis is truncated to the range $[0,90\%]$ for clarity. Comparing these results with those obtained using a larger number of collocated sources and receivers, we observe that increasing the number of measurements improves the stability of the proposed reconstruction method.

\section{Conclusion}\label{sec:conclusion}
In this paper, we investigated Kirchhoff migration imaging for elastic
waves. Using a small-volume expansion, we derived explicit leading-order Kirchhoff migration imaging profiles for the four elastic polarization modes and characterized their resolution. The SS mode provides the highest resolution for localization, whereas the PP mode contains sufficient independent information to recover the EMT and identify the material and geometric parameters of an elliptic inclusion within the leading-order model. For a disk-shaped inclusion, the identification problem is underdetermined, and we established explicit bounds on its volume. We also qualitatively analyzed the error caused by neglected cross-mode contributions, highlighting its dependence on the inclusion location and material parameters. Numerical experiments supported the theoretical predictions and demonstrated the effectiveness of the imaging and reconstruction procedures.

The present analysis assumes a single small inclusion and measurements over a full aperture. Natural extensions include departures from the small-volume regime, incomplete apertures, and multiple inclusions, with further investigation of higher-order scattering, directional resolution, and interactions between inclusions.

\appendix

\section{Far-field asymptotics for $\mGamma^\omega$ and its derivatives}\label{sec:farfield_asymp}
Let
\[
    r:=|\vx|,
    \qquad
    \widehat{\vx}:=\frac{\vx}{|\vx|}.
\]
In the far-field regime $kr\gg1$, the standard large-argument
asymptotic expansion of the Hankel function gives
\[
    H_0^{(1)}(kr)
    \simeq
    \sqrt{\frac{2}{\pi kr}}
    e^{i(kr-\frac{\pi}{4})}.
\]
For the radial function
\[
    f(r):=H_0^{(1)}(kr),
\]
we have
\begin{align*}
    \partial_i f(r)
    &=
    f'(r)\frac{x_i}{r},\\
    \partial_i\partial_j f(r)
    &=
    \frac{f'(r)}{r}\delta_{ij}
    +
    \left(
        f''(r)-\frac{f'(r)}{r}
    \right)
    \frac{x_i x_j}{r^2}.
\end{align*}
Moreover,
\begin{align*}
    \frac{d}{du}H_0^{(1)}(u)
    &=
    -H_1^{(1)}(u),\\
    \frac{d}{du}H_1^{(1)}(u)
    &=
    H_0^{(1)}(u)
    -
    \frac{1}{u}H_1^{(1)}(u).
\end{align*}
Combining these identities with the standard large-argument
asymptotics of the Hankel functions yields
\begin{align}
    \partial_i f(r)
    &\simeq
    ik\widehat{x}_i
    \sqrt{\frac{2}{\pi kr}}
    e^{i(kr-\frac{\pi}{4})},
    \notag\\
    \partial_i\partial_j f(r)
    &\simeq
    -k^2\widehat{x}_i\widehat{x}_j
    \sqrt{\frac{2}{\pi kr}}
    e^{i(kr-\frac{\pi}{4})}.
    \label{eq:hankel_derivative_asymptotics}
\end{align}
Here and below, only the leading terms arising from differentiation
of the oscillatory phase are retained.

For subsequent calculations, we introduce the notation
\[
    \vr_{yz}:=\vy-\vz,
    \qquad
    r_{yz}:=|\vr_{yz}|,
    \qquad
    \widehat{\vr}_{yz}:=\frac{\vr_{yz}}{r_{yz}},
\]
and write
\[
    \widehat{\vr}_{yz}
    =
    \bigl(
        \widehat{r}_{yz,1},
        \widehat{r}_{yz,2}
    \bigr).
\]
Using \eqref{eq:hankel_derivative_asymptotics} in the expression for
the elastic Green tensor, we obtain
\begin{align}\label{eq:gammae_j}
    \mGamma^\omega(\vy-\vz)\ve_j
    \simeq{}&
    -\frac{i}{4\mu}
    \sqrt{\frac{2}{\pi\kappa_S r_{yz}}}
    e^{i(\kappa_Sr_{yz}-\frac{\pi}{4})}
    \left(
        \ve_j
        -
        \widehat{r}_{yz,j}\widehat{\vr}_{yz}
    \right)
    \notag\\
    &-
    \frac{i}{4(\lambda+2\mu)}
    \sqrt{\frac{2}{\pi\kappa_P r_{yz}}}
    e^{i(\kappa_Pr_{yz}-\frac{\pi}{4})}
    \widehat{r}_{yz,j}\widehat{\vr}_{yz}.
\end{align}
Equivalently, its $l$th component satisfies
\begin{align*}
    \left(
        \mGamma^\omega(\vy-\vz)\ve_j
    \right)_l
    \simeq{}&
    -\frac{i}{4\mu}
    \sqrt{\frac{2}{\pi\kappa_S r_{yz}}}
    e^{i(\kappa_Sr_{yz}-\frac{\pi}{4})}
    \left(
        \delta_{jl}
        -
        \widehat{r}_{yz,j}\widehat{r}_{yz,l}
    \right)
    \notag\\
    &-
    \frac{i}{4(\lambda+2\mu)}
    \sqrt{\frac{2}{\pi\kappa_P r_{yz}}}
    e^{i(\kappa_Pr_{yz}-\frac{\pi}{4})}
    \widehat{r}_{yz,j}\widehat{r}_{yz,l}.
\end{align*}

Differentiating the leading-order expression with respect to $z_i$
and retaining the dominant terms gives
\begin{align}\label{eq:dgamma_i}
    \partial_{z_i}
    \left(
        \mGamma^\omega(\vy-\vz)\ve_j
    \right)_l
    \simeq{}&
    -\frac{\kappa_S}{4\mu}
    \sqrt{\frac{2}{\pi\kappa_S r_{yz}}}
    e^{i(\kappa_Sr_{yz}-\frac{\pi}{4})}
    \widehat{r}_{yz,i}
    \left(
        \delta_{jl}
        -
        \widehat{r}_{yz,j}\widehat{r}_{yz,l}
    \right)
    \notag\\
    &-
    \frac{\kappa_P}{4(\lambda+2\mu)}
    \sqrt{\frac{2}{\pi\kappa_P r_{yz}}}
    e^{i(\kappa_Pr_{yz}-\frac{\pi}{4})}
    \widehat{r}_{yz,i}
    \widehat{r}_{yz,j}
    \widehat{r}_{yz,l}.
\end{align}
Similarly,
\begin{align}\label{eq:dgamma_p}
    \partial_{z_p}
    \Gamma_{jq}^\omega(\vy-\vz)
    \simeq{}&
    -\frac{\kappa_S}{4\mu}
    \sqrt{\frac{2}{\pi\kappa_S r_{yz}}}
    e^{i(\kappa_Sr_{yz}-\frac{\pi}{4})}
    \widehat{r}_{yz,p}
    \left(
        \delta_{jq}
        -
        \widehat{r}_{yz,j}\widehat{r}_{yz,q}
    \right)
    \notag\\
    &-
    \frac{\kappa_P}{4(\lambda+2\mu)}
    \sqrt{\frac{2}{\pi\kappa_P r_{yz}}}
    e^{i(\kappa_Pr_{yz}-\frac{\pi}{4})}
    \widehat{r}_{yz,p}
    \widehat{r}_{yz,j}
    \widehat{r}_{yz,q}.
\end{align}

\section{Proofs of results in Section \ref{sec:psf_resolution}}
In this section, we derive the imaging functionals for the different polarization modes. Throughout this section, we adopt the notation introduced in \ref{sec:farfield_asymp}.

\subsection{Proof of Theorem \ref{thm:psf_density}}\label{sec:proof_I^D}
 We calculate the four polarization modes one by one in the following subsections.
\subsubsection{Calculation of $\fI_{P,P}^{\mathrm{sv},\mathrm{D}}$}
We start with the calculation of the integral
\begin{align}
    \mF_{P} = \int_{\partial\Omega} \mGamma^\omega(\vx-\vz)e^{-i\kappa_{P}|\vz^s-\vx|} d\sigma(\vx).
\end{align}
By (\ref{eq:gammae_j}), we have 
\begin{align}
    &\mGamma^\omega(\vx-\vz)e^{-i\kappa_{P}|\vz^s-\vx|} \notag\\
    &\quad\simeq -ie^{-\frac{\pi}{4}i}\bra{\sqrt{\frac{2}{\pi\kappa_{S}r_{xz}}}\frac{\mI-\hat{\vr}_{xz}\hat{\vr}_{xz}^\top}{4\mu}e^{i(\kappa_{S}|\vx-\vz|-\kappa_{P}|\vz^s-\vx|)}
    +
    \sqrt{\frac{2}{\pi\kappa_{P}r_{xz}}}\frac{\hat{\vr}_{xz}\hat{\vr}_{xz}^\top}{4(\lambda+2\mu)}e^{i\kappa_{P}(|\vx-\vz|-|\vz^s-\vx|)}
    }
\end{align}
Then, we observe that in the integral $\mF_{P}$, there are two rapid phase terms:
\begin{align}
    \Phi_1(\vx,\vz) = \frac{\widetilde{\omega}}{\varepsilon^{\frac{1}{2}}}\bra{\frac{|\vx-\vz|}{c_S}-\frac{|\vz^s-\vx|}{c_P}}, \quad \Phi_2(\vx,\vz) = \frac{\widetilde{\omega}}{\varepsilon^{\frac{1}{2}}}\bra{\frac{|\vx-\vz|}{c_P}-\frac{|\vz^s-\vx|}{c_P}}.
\end{align}
By stationary phase analysis, we have
\begin{align}
    \mF_{P} \simeq \frac{-ie^{-\frac{\pi}{4}i}}{4(\lambda+2\mu)}\sqrt{\frac{2}{\pi\kappa_P}} \cdot\int_{\partial\Omega} \frac{\hat{\vr}_{xz}\hat{\vr}_{xz}^\top}{\sqrt{r_{xz}}} e^{i\kappa_{P}(|\vx-\vz|-|\vz^s-\vx|)}  d\sigma(\vx).
\end{align}
To calculate the integral, we adopt the following parametrization of $\vz^s$:
\begin{align}
    \vz^s = \vz +\varepsilon^{\frac{1}{2}}\vr^s.
\end{align}
    
By Taylor's expansion, we have 
\begin{align}
    |\vz^s-\vx| \simeq |\vz-\vx|- \varepsilon^{\frac{1}{2}}\left\langle \hat{\vr}_{xz},\vr^s \right\rangle.
\end{align}
Further, we have the following asymptotic expansions 
\begin{align}
    r_{xz} \simeq |\vx|\bra{1+\fO\bra{\delta}}, \quad \hat{\vr}_{xz}\simeq \frac{\vx}{|\vx|}\bra{1+\fO\bra{\delta}}.
\end{align}
Hence, we have
\begin{align}
    \mF_{P} \simeq \frac{-ie^{-\frac{\pi}{4}i}}{4(\lambda+2\mu)}\sqrt{\frac{2}{\pi\kappa_P}} \cdot\int_{\partial\Omega} \frac{\hat{\vx}\hat{\vx}^\top}{\sqrt{|\vx|}} e^{i\frac{\widetilde{\omega}}{c_P}\left\langle \hat{\vx},\vr^s \right\rangle}  d\sigma(\vx).
\end{align}
We introduce the parametrization that $\vx = (R\cos\theta,R\sin\theta)$ for $\vx\in\partial\Omega$, and we denote 
\begin{align}
    &A_1 = \int_0^{2\pi}  \cos^2\theta e^{i\frac{\widetilde{\omega}}{c_P}(r^s_1\cos\theta+ r^s_2\sin\theta)}d\theta,\\
    &A_2 = \int_0^{2\pi}  \sin\theta\cos\theta e^{i\frac{\widetilde{\omega}}{c_P}(r^s_1\cos\theta+ r^s_2\sin\theta)}d\theta,\\
    &A_3 = \int_0^{2\pi}  \sin^2\theta e^{i\frac{\widetilde{\omega}}{c_P}(r^s_1\cos\theta+ r^s_2\sin\theta)}d\theta.
\end{align}
We have
\begin{align}\label{eq:F_p1}
    \mF_{P} \simeq \frac{-ie^{-\frac{\pi}{4}i}}{4(\lambda+2\mu)}\sqrt{\frac{2R}{\pi\kappa_P}} \cdot
    \bra{\begin{matrix}
        A_1 & A_2\\
        A_2 & A_3
    \end{matrix}}. 
\end{align}

By Lemma \ref{lem:aux}, we have 
\begin{align*}
    &A_1  = \int_0^{2\pi}  \frac{1+\cos2\theta}{2}\cdot e^{i\frac{\widetilde{\omega}}{c_P}(r^s_1\cos\theta+ r^s_2\sin\theta)}d\theta =\pi J_0\bra{\frac{\widetilde{\omega}}{c_P}r^s}-\pi J_2\bra{\frac{\widetilde{\omega}}{c_P}r^s}\cos2\varphi^s, \\
    &A_2 = \int_0^{2\pi}  \frac{1}{2}\sin2\theta \cdot e^{i\frac{\widetilde{\omega}}{c_P}(r^s_1\cos\theta+ r^s_2\sin\theta)}d\theta = -\pi J_2\bra{\frac{\widetilde{\omega}}{c_P}r^s}\sin2\varphi^s,\\
    &A_3 = \int_0^{2\pi}  \frac{1-\cos2\theta}{2}\cdot e^{i\frac{\widetilde{\omega}}{c_P}(r^s_1\cos\theta+ r^s_2\sin\theta)}d\theta =\pi J_0\bra{\frac{\widetilde{\omega}}{c_P}r^s}+\pi J_2\bra{\frac{\widetilde{\omega}}{c_P}r^s}\cos2\varphi^s.
\end{align*}
Therefore, we have 
\begin{align*}
    \fI_{P,P}^{\mathrm{sv},\mathrm{D}}(\vz^s) &\simeq -\varepsilon^2\omega^2(\rho-\rho')|B|\sum_j\ve_j^\top\mF_{P}^2\ve_j \\
    & = -\varepsilon^2\frac{iR\omega^2(\rho-\rho')|B|}{8(\lambda+2\mu)^2\pi\kappa_P}(A_1^2+2A_2^2+A_3^2)\\
    & = -\varepsilon^2\frac{i\pi R\omega(\rho-\rho')|B|c_P}{4(\lambda+2\mu)^2} \bra{J_0^2+J_2^2}\bra{\frac{\widetilde{\omega}}{c_P}r^s}.
\end{align*}

\subsubsection{Calculation of $\fI_{S,S}^{\mathrm{sv},\mathrm{D}}$}
We consider the integral 
\begin{align}
    \mF_s = \int_{\partial\Omega} \mGamma^\omega(\vx-\vz)e^{-i\kappa_{S}|\vz^s-\vx|} d\sigma(\vx).
\end{align}
By a similar rapid phase argument as in the previous proof, we have
\begin{align}
    \mF_{S}
    &\simeq \frac{-ie^{-\frac{\pi}{4}i}}{4\mu}\sqrt{\frac{2}{\pi\kappa_S}} \cdot\int_{\partial\Omega} \frac{\mI-\hat{\vr}_{xz}\hat{\vr}_{xz}^\top}{\sqrt{r_{xz}}} e^{i\kappa_{S}(|\vx-\vz^s|-|\vz^s-\vx|)}  d\sigma(\vx) \\
    &\simeq \frac{-ie^{-\frac{\pi}{4}i}}{4\mu}\sqrt{\frac{2}{\pi\kappa_S}} \cdot\int_{\partial\Omega} \frac{\mI-\hat{\vx}\hat{\vx}^\top}{\sqrt{|\vx|}} e^{i\frac{\widetilde{\omega}}{c_S}\left\langle \hat{\vx},\vr^s \right\rangle}  d\sigma(\vx).
\end{align}
We introduce the parametrization that $\vx = (R\cos\theta,R\sin\theta)$ for $\vx\in\partial\Omega$, and we denote 
\begin{align}
    &A'_1 = \int_0^{2\pi}  (1-\cos^2\theta) e^{i\frac{\widetilde{\omega}}{c_S}(r^s_1\cos\theta+ r^s_2\sin\theta)}d\theta,\\
    &A'_2 = -\int_0^{2\pi}  \sin\theta\cos\theta e^{i\frac{\widetilde{\omega}}{c_S}(r^s_1\cos\theta+ r^s_2\sin\theta)}d\theta,\\
    &A'_3 = \int_0^{2\pi}  (1-\sin^2\theta) e^{i\frac{\widetilde{\omega}}{c_S}(r^s_1\cos\theta+ r^s_2\sin\theta)}d\theta.
\end{align}
We have
\[
    \mF_{S} \simeq \frac{-ie^{-\frac{\pi}{4}i}}{4\mu}\sqrt{\frac{2R}{\pi\kappa_S}} \cdot
    \bra{\begin{matrix}
        A_1' & A_2'\\
        A_2' & A_3'
    \end{matrix}}. 
\]
By Lemma \ref{lem:aux}, we have 
\begin{align*}
    &A_1'  = \int_0^{2\pi}  \frac{1-\cos2\theta}{2}\cdot e^{i\frac{\widetilde{\omega}}{c_S}(r^s_1\cos\theta+ r^s_2\sin\theta)}d\theta =\pi J_0\bra{\frac{\widetilde{\omega}}{c_S}r^s}+\pi J_2\bra{\frac{\widetilde{\omega}}{c_S}r^s}\cos2\varphi^s, \\
    &A_2' = -\int_0^{2\pi}  \frac{1}{2}\sin2\theta \cdot e^{i\frac{\widetilde{\omega}}{c_S}(r^s_1\cos\theta+ r^s_2\sin\theta)}d\theta = \pi J_2\bra{\frac{\widetilde{\omega}}{c_S}r^s}\sin2\varphi^s,\\
    &A_3' = \int_0^{2\pi}  \frac{1+\cos2\theta}{2}\cdot e^{i\frac{\widetilde{\omega}}{c_S}(r^s_1\cos\theta+ r^s_2\sin\theta)}d\theta =\pi J_0\bra{\frac{\widetilde{\omega}}{c_S}r^s}-\pi J_2\bra{\frac{\widetilde{\omega}}{c_S}r^s}\cos2\varphi^s.
\end{align*}
Therefore, we have 
\begin{align*}
    \fI_{S,S}^{\mathrm{sv},\mathrm{D}}(\vz^s) &\simeq -\varepsilon^2\omega^2(\rho-\rho')|B|\sum_j\ve_j^\top\mF_{S}^2\ve_j \\
    & = -\varepsilon^2\frac{i\pi R\omega(\rho-\rho')|B|c_S}{4\mu^2} \bra{J_0^2+J_2^2}\bra{\frac{\widetilde{\omega}}{c_S}r^s}
\end{align*}
\subsubsection{Calculation of $\fI_{P,S}^{\mathrm{sv},\mathrm{D}}$ and $\fI_{S,P}^{\mathrm{sv},\mathrm{D}}$}
By a similar argument in the previous two proofs, we have 
\begin{align*}
    \fI_{P,S}^{\mathrm{sv},\mathrm{D}}(\vz^s) &\simeq -\varepsilon^2\omega^2(\rho-\rho')|B|\sum_j\ve_j^\top\mF_{P}\mF_{S}\ve_j \\
    & = -\varepsilon^2\frac{i R\omega(\rho-\rho')|B|\sqrt{c_Sc_P}}{8\mu(\lambda+2\mu)\pi}\bra{A_1A_1'+2A_2A_2'+A_3A_3'}\\
    & = -\varepsilon^2\frac{i \pi R\omega(\rho-\rho')|B|\sqrt{c_Sc_P}}{4\mu(\lambda+2\mu)}\bra{J_0\bra{\frac{\widetilde{\omega}}{c_P}r^s}J_0\bra{\frac{\widetilde{\omega}}{c_S}r^s}-J_2\bra{\frac{\widetilde{\omega}}{c_P}r^s}J_2\bra{\frac{\widetilde{\omega}}{c_S}r^s}}.
\end{align*}
Also, we can calculate that
\begin{align*}
    \fI_{P,S}^{\mathrm{sv},\mathrm{D}}(\vz^s) &\simeq -\varepsilon^2\omega^2(\rho-\rho')|B|\sum_j\ve_j^\top \mF_{P}\mF_{S}\ve_j \\
     &= -\varepsilon^2\omega^2(\rho-\rho')|B|\cdot \operatorname{Tr}(\mF_{P}\mF_{S}) \\
     &= -\varepsilon^2\omega^2(\rho-\rho')|B|\cdot \operatorname{Tr}(\mF_{S}\mF_{P}) \\
     &= \fI_{S,P}^{\mathrm{sv},\mathrm{D}}(\vz^s).
\end{align*}

\subsection{Proof of Theorem \ref{thm:psf_elas}}\label{sec:proof_I^E}
We calculate the four polarization modes one by one in the following subsections.
\subsubsection{Calculation of $\fI_{P,P}^{\mathrm{sv},\mathrm{E}}$}
By (\ref{eq:dgamma_i}), (\ref{eq:dgamma_p}), and the rapid phase argument, we have 
\begin{align*}
    \fI_{P,P}^{\mathrm{sv},\mathrm{E}}(\vz^s) &= -\varepsilon^2\sum_{j=1,2}\int_{\partial\Omega} \int_{\partial\Omega} 
    d\sigma(\vx)d\sigma(\vy)~e^{-i\kappa_P(|\vx-\vz^s|+|\vy-\vz^s|)} \\
    &\quad\times\sum_{i,l=1,2} \bra{ \partial_{z_i}(\mGamma^\omega(\vz-\vy)\ve_j)_l\sum_{p,q=1,2}m_{ilpq}\partial_{z_p}\Gamma^\omega_{jq}(\vx-\vz)  }\\
    &\simeq -\varepsilon^2\sum_{j=1,2}\frac{\kappa_P^2}{16(\lambda+2\mu)^2} \int_{\partial\Omega} \int_{\partial\Omega} 
    d\sigma(\vx)d\sigma(\vy)~\frac{-2i}{\pi\kappa_P \sqrt{r_{yz}r_{xz}}}e^{-i\kappa_P(|\vx-\vz^s|+|\vy-\vz^s|)}e^{i\kappa_P(r_{yz}+r_{xz})}\\
    &\quad\times\sum_{i,l}\hat{r}_{yz,i}\hat{r}_{yz,j}\hat{r}_{yz,l}\sum_{p,q}m_{ilpq}\hat{r}_{xz,p}\hat{r}_{xz,j}\hat{r}_{xz,q} \\
    &\simeq \frac{i\kappa_PR\varepsilon^2}{8\pi(\lambda+2\mu)^2}\sum_{j,i,l,p,q}m_{ilpq}\bra{\int_{\partial\Omega} \frac{\hat{y}_i\hat{y}_j\hat{y}_l}{\sqrt{|\vy|}} e^{i\frac{\widetilde{\omega}}{c_P}\left\langle \hat{\vy},\vr^s \right\rangle} d\sigma(\vy)}
    \bra{\int_{\partial\Omega} \frac{\hat{x}_p\hat{x}_j\hat{x}_q}{\sqrt{|\vx|}} e^{i\frac{\widetilde{\omega}}{c_P}\left\langle \hat{\vx},\vr^s \right\rangle} d\sigma(\vx)}
\end{align*}
We calculate that
\begin{align*}
    &\int_0^{2\pi}  \cos^3\theta 
    e^{i\frac{\widetilde{\omega}}{c_P}(r^s_1\cos\theta+ r^s_2\sin\theta)}d\theta
    = \frac{i\pi}{2}\Bigl(
    3J_1\!\left(\frac{\widetilde{\omega}}{c_P}r^s\right)\cos\varphi^s
    -J_3\!\left(\frac{\widetilde{\omega}}{c_P}r^s\right)\cos3\varphi^s
    \Bigr),\\
    &\int_0^{2\pi}  \sin\theta \cos^2\theta 
    e^{i\frac{\widetilde{\omega}}{c_P}(r^s_1\cos\theta+ r^s_2\sin\theta)}d\theta
    = \frac{i\pi}{2}\Bigl(
    J_1\!\left(\frac{\widetilde{\omega}}{c_P}r^s\right)\sin\varphi^s
    -J_3\!\left(\frac{\widetilde{\omega}}{c_P}r^s\right)\sin3\varphi^s
    \Bigr),\\
    &\int_0^{2\pi}  \sin^2\theta\cos\theta 
    e^{i\frac{\widetilde{\omega}}{c_P}(r^s_1\cos\theta+ r^s_2\sin\theta)}d\theta
    = \frac{i\pi}{2}\Bigl(
    J_1\!\left(\frac{\widetilde{\omega}}{c_P}r^s\right)\cos\varphi^s
    +J_3\!\left(\frac{\widetilde{\omega}}{c_P}r^s\right)\cos3\varphi^s
    \Bigr),\\ 
    &\int_0^{2\pi}  \sin^3\theta 
    e^{i\frac{\widetilde{\omega}}{c_P}(r^s_1\cos\theta+ r^s_2\sin\theta)}d\theta
    = \frac{i\pi}{2}\Bigl(
    3J_1\!\left(\frac{\widetilde{\omega}}{c_P}r^s\right)\sin\varphi^s
    +J_3\!\left(\frac{\widetilde{\omega}}{c_P}r^s\right)\sin3\varphi^s
    \Bigr),
\end{align*}
where we use the notation that $\varphi^s = \operatorname{arg}(r_1^s+ir^s_2)$.
Therefore, we derive
\begin{align*}
   \fI_{P,P}^{\mathrm{sv},\mathrm{E}}(\vz^s) \simeq\frac{-i\kappa_P\pi R\varepsilon^2}{32(\lambda+2\mu)^2}&\left[M_{p,p}^{(1)}J_1^2\left(\frac{\widetilde{\omega}}{c_P}r^s\right)+M_{p,p}^{(2)}J_3^2\left(\frac{\widetilde{\omega}}{c_P}r^s\right)+M_{p,p}^{(3)}J_1\left(\frac{\widetilde{\omega}}{c_P}r^s\right)J_3\left(\frac{\widetilde{\omega}}{c_P}r^s\right)\cos4\varphi^s\right. \\
    &~+M_{p,p}^{(4)}J_1^2\left(\frac{\widetilde{\omega}}{c_P}r^s\right)\sin^2\varphi^s+M_{p,p}^{(5)}J_1^2\left(\frac{\widetilde{\omega}}{c_P}r^s\right)\cos^2\varphi^s\\
    &~\left.+M_{p,p}^{(6)}J_1\left(\frac{\widetilde{\omega}}{c_P}r^s\right)J_3\left(\frac{\widetilde{\omega}}{c_P}r^s\right)\cos\varphi^s\cos 3\varphi^s+M_{p,p}^{(7)}J_1\left(\frac{\widetilde{\omega}}{c_P}r^s\right)J_3\left(\frac{\widetilde{\omega}}{c_P}r^s\right)\sin\varphi^s\sin3\varphi^s\right]
\end{align*}
where $M_{p,p}^{(1)} = 4m_{1212}+6m_{1122}$, $M_{p,p}^{(2)} = m_{1111}+m_{2222}+4m_{1212}-2m_{1122}$, $M_{p,p}^{(3)} = 8m_{1212}+4m_{1122}$, $M_{p,p}^{(4)}=9m_{2222}+m_{1111}$, $M_{p,p}^{(5)}=9m_{1111}+m_{2222}$, $M_{p,p}^{(6)}=-6m_{1111}+2m_{2222}$, $M_{p,p}^{(7)}=6m_{2222}-2m_{1111}$.
Hence, we derive that 
\begin{align*}
    \fI_{P,P}^{\mathrm{sv},\mathrm{E}}(\vz^s) \simeq\frac{-i\kappa_P\pi R\varepsilon^2}{32(\lambda+2\mu)^2}\bra{A_{P,P}^{(1)}(\varphi^s)
     J_1^2\!\left(\frac{\widetilde{\omega}}{c_S}r^s\right) 
    + A_{P,P}^{(2)}(\varphi^s) J_3^2\!\left(\frac{\widetilde{\omega}}{c_S}r^s\right) 
    +A_{P,P}^{(3)}(\varphi^s)J_1\!\left(\frac{\widetilde{\omega}}{c_S}r^s\right) J_3\!\left(\frac{\widetilde{\omega}}{c_S}r^s\right) },
\end{align*}
where the coefficients are given by 
\begin{align*}
    &A_{P,P}^{(1)}(\varphi^s) = 4m_{1212}+6m_{1122}+(9m_{2222}+m_{1111})\sin^2\varphi^s+(9m_{1111}+m_{2222})\cos^2\varphi^s, \\
    &A_{P,P}^{(2)}(\varphi^s) = m_{1111}+m_{2222}+4m_{1212}-2m_{1122}, \\
    &A_{P,P}^{(3)}(\varphi^s) = (8m_{1212}+4m_{1122})\cos4\varphi^s-(6m_{1111}-2m_{2222})\cos\varphi^s\cos3\varphi^s+(6m_{2222}-2m_{1111})\sin\varphi^s\sin3\varphi^s.
\end{align*}
\subsubsection{Calculation of $\fI_{S,S}^{\mathrm{sv},\mathrm{E}}$}
By the rapid phase argument, we have 
\begin{align*}
    \fI_{S,S}^{\mathrm{sv},\mathrm{E}}(\vz^s) &= -\varepsilon^2\sum_{j=1,2}\int_{\partial\Omega} \int_{\partial\Omega} 
    d\sigma(\vx)d\sigma(\vy)~e^{-i\kappa_S(|\vx-\vz^s|+|\vy-\vz^s|)} \\
    &\quad\times\sum_{i,l=1,2} \bra{ \partial_{z_i}(\mGamma^\omega(\vz-\vy)\ve_j)_l\sum_{p,q=1,2}m_{ilpq}\partial_{z_p}\Gamma^\omega_{jq}(\vx-\vz)  }\\
    &\simeq -\varepsilon^2\sum_{j=1,2}\frac{\kappa_S^2}{16\mu^2} \int_{\partial\Omega} \int_{\partial\Omega} 
    d\sigma(\vx)d\sigma(\vy)~\frac{-2i}{\pi\kappa_S \sqrt{r_{yz}r_{xz}}}e^{-i\kappa_S(|\vx-\vz^s|+|\vy-\vz^s|)}e^{i\kappa_S(r_{yz}+r_{xz})}\\
    &\quad\times\sum_{i,l}\hat{r}_{yz,i}\bra{\delta_{jl}-\hat{r}_{yz,j}\hat{r}_{yz,l}}\sum_{p,q}m_{ilpq}\hat{r}_{xz,p}(\delta_{jq}-\hat{r}_{xz,j}\hat{r}_{xz,q}) \\
\end{align*}
We calculate that
\begin{align*}
    &\sum_{i,l}\hat{r}_{yz,i}\bra{\delta_{jl}-\hat{r}_{yz,j}\hat{r}_{yz,l}}\sum_{p,q}m_{ilpq}\hat{r}_{xz,p}(\delta_{jq}-\hat{r}_{xz,j}\hat{r}_{xz,q}) 
    = 
    \sum_i \hat{r}_{yz,i}\sum_pm_{ijpj}\hat{r}_{xz,p}\\
    &\quad - \sum_{i,l}\hat{r}_{yz,i}\hat{r}_{yz,j}\hat{r}_{yz,l}\sum_{p}m_{ilpj}\hat{r}_{xz,p}-\sum_{i}\hat{r}_{yz,i}\sum_{p,q}m_{ijpq}\hat{r}_{xz,p}\hat{r}_{xz,j}\hat{r}_{xz,q}\\
    &\quad + \sum_{i,l}\hat{r}_{yz,i}\hat{r}_{yz,j}\hat{r}_{yz,l}\sum_{p,q}m_{ilpq}\hat{r}_{xz,p}\hat{r}_{xz,j}\hat{r}_{xz,q}.
\end{align*}
Adopting the asymptotic regime, we can calculate the imaging functional as
\begin{align*}
    \fI_{S,S}^{\mathrm{sv},\mathrm{E}}(\vz^s) &\simeq \frac{i\kappa_SR\varepsilon^2}{8\pi\mu^2}\cdot\left(-4\pi^2(m_{1111}\cos^2\varphi^s+m_{2222}\sin^2\varphi^s+m_{1212})J_1^2\bra{\frac{\widetilde{\omega}}{c_S}r^s} \right.\\
    &\quad+2\pi^2(3m_{1111}\cos^2\varphi^s+3m_{2222}\sin^2\varphi^s+m_{1122}+2m_{1212})J_1^2\bra{\frac{\widetilde{\omega}}{c_S}r^s}\\
    &\quad+2\pi^2(m_{1122}+2m_{1212})J_1\bra{\frac{\widetilde{\omega}}{c_S}r^s}J_3\bra{\frac{\widetilde{\omega}}{c_S}r^s}\cos4\varphi^s\\
    &\quad+2\pi^2(-m_{1111}\cos\varphi^s\cos3\varphi^s+m_{2222}\sin\varphi^s\sin3\varphi^s )J_1\bra{\frac{\widetilde{\omega}}{c_S}r^s}J_3\bra{\frac{\widetilde{\omega}}{c_S}r^s} \\
    &\quad-\frac{\pi^2}{4}\left[M_{p,p}^{(1)}J_1^2\left(\frac{\widetilde{\omega}}{c_S}r^s\right)+M_{p,p}^{(2)}J_3^2\left(\frac{\widetilde{\omega}}{c_S}r^s\right)+M_{p,p}^{(3)}J_1\left(\frac{\widetilde{\omega}}{c_S}r^s\right)J_3\left(\frac{\widetilde{\omega}}{c_S}r^s\right)\cos4\varphi^s\right. \\
    &\quad+M_{p,p}^{(4)}J_1^2\left(\frac{\widetilde{\omega}}{c_S}r^s\right)\sin^2\varphi^s+M_{p,p}^{(5)}J_1^2\left(\frac{\widetilde{\omega}}{c_S}r^s\right)\cos^2\varphi^s\\
    &\left.\quad\left.+M_{p,p}^{(6)}J_1\left(\frac{\widetilde{\omega}}{c_S}r^s\right)J_3\left(\frac{\widetilde{\omega}}{c_S}r^s\right)\cos\varphi^s\cos 3\varphi^s+M_{p,p}^{(7)}J_1\left(\frac{\widetilde{\omega}}{c_S}r^s\right)J_3\left(\frac{\widetilde{\omega}}{c_S}r^s\right)\sin\varphi^s\sin3\varphi^s\right]\right) \\
    &\simeq -\frac{i\kappa_S R\varepsilon^2 \pi}{8\mu^2}
    \bra{A_{S,S}^{(1)}(\varphi^s)
     J_1^2\!\left(\frac{\widetilde{\omega}}{c_S}r^s\right) 
    + A_{S,S}^{(2)}(\varphi^s) J_3^2\!\left(\frac{\widetilde{\omega}}{c_S}r^s\right) 
    +A_{S,S}^{(3)}(\varphi^s)J_1\!\left(\frac{\widetilde{\omega}}{c_S}r^s\right) J_3\!\left(\frac{\widetilde{\omega}}{c_S}r^s\right) },
\end{align*}
where we use the convention that $\varphi^s = \operatorname{arg}(r_1^s+ir^s_2)$, and the coefficients are given by 
\begin{align*}
    & A_{S,S}^{(1)}(\varphi^s) = A_{S,S}^{(2)}(\varphi^s) = \frac{1}{4}m_{1111}+\frac{1}{4}m_{2222}+m_{1212}-\frac{1}{2}m_{1122},\\
    & A_{S,S}^{(3)}(\varphi^s) = \bra{\frac{1}{2}(m_{1111}+m_{2222})-m_{1122}-2m_{1212}}\cos4\varphi^s.
\end{align*}
\subsubsection{Calculation of $\fI_{S,P}^{\mathrm{sv},\mathrm{E}}$ and $\fI_{P,S}^{\mathrm{sv},\mathrm{E}}$}
By the rapid phase argument, we have 
\begin{align*}
    \fI_{S,P}^{\mathrm{sv},\mathrm{E}}(\vz^s) &= -\varepsilon^2\sum_{j=1,2}\int_{\partial\Omega} \int_{\partial\Omega} 
    d\sigma(\vx)d\sigma(\vy)~e^{-i\kappa_S|\vx-\vz^s|-i\kappa_P|\vy-\vz^s|} \\
    &\quad\times\sum_{i,l=1,2} \bra{ \partial_{z_i}(\mGamma^\omega(\vz-\vy)\ve_j)_l\sum_{p,q=1,2}m_{ilpq}\partial_{z_p}\Gamma^\omega_{jq}(\vx-\vz)  }\\
    &\simeq -\varepsilon^2\sum_{j=1,2}\frac{\kappa_S\kappa_P}{16\mu(\lambda+2\mu)} \int_{\partial\Omega} \int_{\partial\Omega} 
    d\sigma(\vx)d\sigma(\vy)~\frac{-2i}{\pi\sqrt{\kappa_P\kappa_S} \sqrt{r_{yz}r_{xz}}}e^{-i\kappa_S(|\vx-\vz^s|-r_{xz})}e^{-i\kappa_P(|\vy-\vz^s|-r_{yz})}\\
    &\quad\times\sum_{i,l}\hat{r}_{yz,i}\hat{r}_{yz,j}\hat{r}_{yz,l}\sum_{p,q}m_{ilpq}\hat{r}_{xz,p}(\delta_{jq}-\hat{r}_{xz,j}\hat{r}_{xz,q}) \\
\end{align*}
We calculate that
\begin{align*}
    &\sum_{i,l}\hat{r}_{yz,i}\hat{r}_{yz,j}\hat{r}_{yz,l}\sum_{p,q}m_{ilpq}\hat{r}_{xz,p}(\delta_{jq}-\hat{r}_{xz,j}\hat{r}_{xz,q})\\
    &\quad=  \sum_{i,l}\hat{r}_{yz,i}\hat{r}_{yz,j}\hat{r}_{yz,l}\sum_{p}m_{ilpj}\hat{r}_{xz,p}- \sum_{i,l}\hat{r}_{yz,i}\hat{r}_{yz,j}\hat{r}_{yz,l}\sum_{p,q}m_{ilpq}\hat{r}_{xz,p}\hat{r}_{xz,j}\hat{r}_{xz,q}.
\end{align*}
Adopting the asymptotic regime, we can calculate the imaging functional as
\begin{align*}
    \fI_{S,P}^{\mathrm{sv},\mathrm{E}}(\vz^s) &\simeq \frac{i\sqrt{\kappa_S\kappa_P}R\varepsilon^2}{8\pi\mu(\lambda{+2\mu})}\cdot\left(
    -\pi^2(3m_{1111}\cos^2\varphi^s+3m_{2222}\sin^2\varphi^s+m_{1122}+2m_{1212})J_1\bra{\frac{\widetilde{\omega}}{c_S}r^s}J_1\bra{\frac{\widetilde{\omega}}{c_P}r^s}\right.\\
    &\quad-\pi^2(m_{1122}+2m_{1212})J_3\bra{\frac{\widetilde{\omega}}{c_P}r^s}J_1\bra{\frac{\widetilde{\omega}}{c_S}r^s}\cos4\varphi^s\\
    &\quad-\pi^2(-m_{1111}\cos\varphi^s\cos3\varphi^s+m_{2222}\sin\varphi^s\sin3\varphi^s )J_3\bra{\frac{\widetilde{\omega}}{c_P}r^s}J_1\bra{\frac{\widetilde{\omega}}{c_S}r^s} \\
    &\quad+\frac{\pi^2}{4}\left[M_{p,p}^{(1)}J_1\left(\frac{\widetilde{\omega}}{c_S}r^s\right)J_1\left(\frac{\widetilde{\omega}}{c_P}r^s\right)+M_{p,p}^{(2)}J_3\left(\frac{\widetilde{\omega}}{c_S}r^s\right)J_3\left(\frac{\widetilde{\omega}}{c_P}r^s\right)\right.\\
    &\quad+\frac{M_{p,p}^{(3)}}{2}J_1\left(\frac{\widetilde{\omega}}{c_S}r^s\right)J_3\left(\frac{\widetilde{\omega}}{c_P}r^s\right)\cos4\varphi^s+\frac{M_{p,p}^{(3)}}{2}J_1\left(\frac{\widetilde{\omega}}{c_P}r^s\right)J_3\left(\frac{\widetilde{\omega}}{c_S}r^s\right)\cos4\varphi^s \\
    &\quad+M_{p,p}^{(4)}J_1\left(\frac{\widetilde{\omega}}{c_S}r^s\right)J_1\left(\frac{\widetilde{\omega}}{c_P}r^s\right)\sin^2\varphi^s+M_{p,p}^{(5)}J_1\left(\frac{\widetilde{\omega}}{c_S}r^s\right)J_1\left(\frac{\widetilde{\omega}}{c_P}r^s\right)\cos^2\varphi^s\\
    &\quad+\frac{M_{p,p}^{(6)}}{2}J_1\left(\frac{\widetilde{\omega}}{c_S}r^s\right)J_3\left(\frac{\widetilde{\omega}}{c_P}r^s\right)\cos\varphi^s\cos 3\varphi^s+\frac{M_{p,p}^{(6)}}{2}J_1\left(\frac{\widetilde{\omega}}{c_P}r^s\right)J_3\left(\frac{\widetilde{\omega}}{c_S}r^s\right)\cos\varphi^s\cos 3\varphi^s\\
    &\quad\left.\left.+\frac{M_{p,p}^{(7)}}{2}J_1\left(\frac{\widetilde{\omega}}{c_S}r^s\right)J_3\left(\frac{\widetilde{\omega}}{c_P}r^s\right)\sin\varphi^s\sin3\varphi^s +\frac{M_{p,p}^{(7)}}{2}J_1\left(\frac{\widetilde{\omega}}{c_P}r^s\right)J_3\left(\frac{\widetilde{\omega}}{c_S}r^s\right)\sin\varphi^s\sin 3\varphi^s\right]\right) \\
    &=  -\frac{i\sqrt{\kappa_S\kappa_P}R\varepsilon^2\pi}{8\mu(\lambda+2\mu)}
    \left(A_{S,P}^{(1)}(\varphi^s)
     J_1\!\left(\frac{\widetilde{\omega}}{c_S}r^s\right) J_1\!\left(\frac{\widetilde{\omega}}{c_P}r^s\right) 
    + A_{S,P}^{(2)}(\varphi^s) J_3\!\left(\frac{\widetilde{\omega}}{c_S}r^s\right) J_3\!\left(\frac{\widetilde{\omega}}{c_P}r^s\right) \right.\\
    &\left.\quad\quad+A_{S,P}^{(3)}(\varphi^s)J_1\!\left(\frac{\widetilde{\omega}}{c_S}r^s\right) J_3\!\left(\frac{\widetilde{\omega}}{c_P}r^s\right) + A_{S,P}^{(4)}(\varphi^s)J_1\!\left(\frac{\widetilde{\omega}}{c_P}r^s\right) J_3\!\left(\frac{\widetilde{\omega}}{c_S}r^s\right) \right),
\end{align*}
where the coefficients are given by 
\begin{align*}
&A_{S,P}^{(1)}(\varphi^s) = \frac{1}{4}(m_{1111}+m_{2222}) + m_{1212} - \frac{1}{2}m_{1122} + \frac{1}{2}(m_{1111}-m_{2222})\cos2\varphi^s,\\
&A_{S,P}^{(2)}(\varphi^s) = -\frac{1}{4}(m_{1111}+m_{2222}) - m_{1212} + \frac{1}{2}m_{1122},\\
&A_{S,P}^{(3)}(\varphi^s) = \left( m_{1212} + \frac{1}{2}m_{1122} - \frac{1}{4}(m_{1111}+m_{2222}) \right)\cos4\varphi^s,\\
&A_{S,P}^{(4)}(\varphi^s) = \left( \frac{1}{4}(m_{1111}+m_{2222}) - m_{1212} - \frac{1}{2}m_{1122} \right)\cos4\varphi^s + \frac{1}{2}(m_{1111}-m_{2222})\cos2\varphi^s.
\end{align*}
Interchanging the integration variables \(\vx\) and \(\vy\) and using the major symmetry \(m_{ilpq}=m_{pqil}\) of the EMT, we obtain
\[
I_{S,P}^{\mathrm{sv},E}(\vz^s)
=
I_{P,S}^{\mathrm{sv},E}(\vz^s).
\]
When the inclusion is disk-shaped, the coefficients $\{A_{\alpha,\beta}^{(j)}(\varphi^s)\}_{j,\alpha,\beta}$ can be simplified using (\ref{eq:m_1})-(\ref{eq:m_2}).

\section{Proof of results in Section \ref{sec:identification}}
\subsection{Proof of Theorem~\ref{thm:emt_uniqueness}}
\label{proof:emt_uniqueness}

For simplicity of notation, we introduce the nonzero constants
\[
    C_0
    :=
    -\frac{i\pi R\omega c_P\varepsilon^2}
           {4(\lambda+2\mu)^2},
    \qquad
    C_1
    :=
    -\frac{i\kappa_P\pi R\varepsilon^2}
           {8(\lambda+2\mu)^2},
\]
where $R$ denotes the radius of the measurement curve. We also define
\begin{align*}
    h_2(R_s)
    &:=
    J_1^2\left(\frac{\omega R_s}{c_P}\right)
    -
    J_1\left(\frac{\omega R_s}{c_P}\right)
    J_3\left(\frac{\omega R_s}{c_P}\right),\\
    h_4(R_s)
    &:=
    J_1\left(\frac{\omega R_s}{c_P}\right)
    J_3\left(\frac{\omega R_s}{c_P}\right).
\end{align*}

For a fixed sampling radius $R_s>0$, a sampling point in the
laboratory coordinate system is parameterized as
\[
    \vz^s
    =
    \vz
    +
    R_s
    \begin{pmatrix}
        \cos\varphi^s\\
        \sin\varphi^s
    \end{pmatrix},
\]
where $\varphi^s$ denotes the sampling angle measured from the
$x_1$-axis of the laboratory coordinate system. The corresponding
sampling angle in the principal-axis coordinate system is
\[
    \theta^s
    :=
    \varphi^s-\phi
    \quad(\operatorname{mod}2\pi).
\]
Consequently, the leading-order PP imaging functional
\eqref{eq:both_contrast} can be written as
\begin{equation}\label{eq:angular_decomposition}
\begin{aligned}
    &\fI^{\mathrm{asym}}_{P,P}
    \left(
        \vz+
        R_s
        \begin{pmatrix}
            \cos\varphi^s\\
            \sin\varphi^s
        \end{pmatrix}
    \right)
    =
    A_0(R_s)
    +
    A_2(R_s)\cos(2\theta^s)
    +
    A_4(R_s)\cos(4\theta^s),
\end{aligned}
\end{equation}
where
\begin{align}
    A_0(R_s)
    &=
    C_0M^{(0)}g_0(R_s)
    +
    C_1\bigl[
        M^{(1)}g_1(R_s)
        +M^{(3)}g_3(R_s)
    \bigr],
    \label{eq:A0_coeff}\\
    A_2(R_s)
    &=
    C_1M^{(2)}h_2(R_s),
    \label{eq:A2_coeff}\\
    A_4(R_s)
    &=
    C_1M^{(4)}h_4(R_s).
    \label{eq:A4_coeff}
\end{align}
For each laboratory-frame sampling angle $\varphi_n^s$, recall that
\[
    \theta_n^s
    :=
    \varphi_n^s-\phi
    \quad(\operatorname{mod}2\pi).
\]
At each sampling radius $R_m$, evaluating
\eqref{eq:angular_decomposition} at the angles
$\{\varphi_n^s\}_{n=1}^{N_{\varphi^s}}$ gives
\[
    \begin{pmatrix}
      \fI^{\mathrm{asym}}_{P,P}
      \left(
        \vz+
        R_m
        \begin{pmatrix}
          \cos\varphi_1^s\\
          \sin\varphi_1^s
        \end{pmatrix}
      \right)\\
      \vdots\\
      \fI^{\mathrm{asym}}_{P,P}
      \left(
        \vz+
        R_m
        \begin{pmatrix}
          \cos\varphi_{N_{\varphi^s}}^s\\
          \sin\varphi_{N_{\varphi^s}}^s
        \end{pmatrix}
      \right)
    \end{pmatrix}
    =
    \mV
    \begin{pmatrix}
      A_0(R_m)\\
      A_2(R_m)\\
      A_4(R_m)
    \end{pmatrix},
\]
where
\[
    \mV
    =
    \begin{pmatrix}
      1 & \cos(2\theta_1^s) & \cos(4\theta_1^s)\\
      \vdots & \vdots & \vdots\\
      1 & \cos(2\theta_{N_{\varphi^s}}^s)
        & \cos(4\theta_{N_{\varphi^s}}^s)
    \end{pmatrix}.
\]
Since $\mV$ has rank three, the angular coefficients
$A_0(R_m)$, $A_2(R_m)$, and $A_4(R_m)$ are uniquely determined from
the angular samples.

By assumption, there exists a sampling radius $R_{m_2}$ such that
$h_2(R_{m_2})\neq0$. It follows from \eqref{eq:A2_coeff} that
\[
    M^{(2)}
    =
    \frac{A_2(R_{m_2})}
         {C_1h_2(R_{m_2})}.
\]
Similarly, since $h_4(R_{m_4})\neq0$ for some sampling radius
$R_{m_4}$, equation~\eqref{eq:A4_coeff} yields
\[
    M^{(4)}
    =
    \frac{A_4(R_{m_4})}
         {C_1h_4(R_{m_4})}.
\]

It remains to determine $M^{(0)}$, $M^{(1)}$, and $M^{(3)}$.
Evaluating \eqref{eq:A0_coeff} at the three sampling radii
$R_1,R_2,R_3$ gives
\[
    \begin{pmatrix}
      A_0(R_1)\\
      A_0(R_2)\\
      A_0(R_3)
    \end{pmatrix}
    =
    \begin{pmatrix}
      C_0g_0(R_1) & C_1g_1(R_1) & C_1g_3(R_1)\\
      C_0g_0(R_2) & C_1g_1(R_2) & C_1g_3(R_2)\\
      C_0g_0(R_3) & C_1g_1(R_3) & C_1g_3(R_3)
    \end{pmatrix}
    \begin{pmatrix}
      M^{(0)}\\
      M^{(1)}\\
      M^{(3)}
    \end{pmatrix}.
\]
The determinant of the coefficient matrix is
\[
    C_0C_1^2
    \det
    \begin{pmatrix}
      g_0(R_1) & g_1(R_1) & g_3(R_1)\\
      g_0(R_2) & g_1(R_2) & g_3(R_2)\\
      g_0(R_3) & g_1(R_3) & g_3(R_3)
    \end{pmatrix},
\]
which is nonzero by the radial sampling assumption. Hence,
$M^{(0)}$, $M^{(1)}$, and $M^{(3)}$ are uniquely determined.
Together with the preceding recovery of $M^{(2)}$ and $M^{(4)}$,
this uniquely determines all five coefficients
$M^{(0)},M^{(1)},\ldots,M^{(4)}$.

Moreover, the four EMT-dependent coefficients satisfy
\begin{align}\label{eq:coeff_to_EMT}
\begin{pmatrix}
M^{(1)}\\
M^{(2)}\\
M^{(3)}\\
M^{(4)}
\end{pmatrix}
&=
\begin{pmatrix}
\frac54&\frac54&\frac32&1\\
1&-1&0&0\\
\frac14&\frac14&-\frac12&1\\
-\frac12&-\frac12&1&2
\end{pmatrix}
\begin{pmatrix}
m_{1111}\\
m_{2222}\\
m_{1122}\\
m_{1212}
\end{pmatrix} \triangleq
\mM_c
\begin{pmatrix}
m_{1111}\\
m_{2222}\\
m_{1122}\\
m_{1212}
\end{pmatrix}.
\end{align}
The coefficient matrix $\mM_c$ is invertible. Therefore,
$M^{(1)},\ldots,M^{(4)}$ uniquely determine the EMT expressed in the
principal-axis coordinate system. This completes the proof.

\subsection{Proof of Theorem \ref{thm:volume}}\label{proof:elliptic_thm}
The three constraints given in (\ref{eq:constraint_lame_inclu}) play an important role in the derivation of the estimates of the inclusion volume.\\
\textbf{Step 1.}
Let 
\[
    u = \mu'-\mu, \quad v = \lambda'-\lambda,
\]
denote 
\[
    A = \lambda+2\mu, \quad C = \mu(\tau+1), \quad D = \lambda+\mu,
\]
and denote 
\[  
    m_1 = m_{1122}, \quad m_2= m_{1212}
\]
for simplicity. Then, we have
\[
    m_1 = |B|\frac{A(v(C+u)+u^2)}{(A+u+v)(C+\tau u)}, \quad m_2 = |B| \frac{Cu}{C+\tau u}.
\]
We can then solve $u$ from $m_2$ as:
\begin{align}
    u = \frac{m_2C}{\Delta}, \quad\Delta=|B|C-m_2\tau.
\end{align}

We express $v$ using $u$ through the expression $m_1$ as:
\begin{align}
        v = -\frac{u(ku-L)}{ku+m_2AC},
\end{align}
where
\[  
    k = m_2A-m_1C, \quad L = m_1AC.
\]
\textbf{Step 2.} Constraint $\mu'>0$:\\
By the calculation in the previous step, We have 
\[
    u = \frac{m_2C}{\Delta}>-\mu.
\]
When $\Delta>0$ $(|B|>\frac{m_2\tau}{C})$, we have
\[
    m_2C>-\mu\Delta=-\mu(|B|C-m_2\tau) \Rightarrow |B|>-\frac{m_2}{C}.
\]
Thus, we have 
\[
    |B|>\max \bra{\frac{m_2\tau}{C},-\frac{m_2}{C}}.
\]
When $\Delta<0$ $(|B|<\frac{m_2\tau}{C})$, we have
\[
    |B|<\min \bra{\frac{m_2\tau}{C},-\frac{m_2}{C}}.
\]
We further notice that $|B|\ge 0$, the constraint $\mu'>0$ implies the following constraint on $|B|$:
\begin{align}
    |B|>\max \bra{\frac{m_2\tau}{\mu(\tau+1)},-\frac{m_2}{\mu(\tau+1)}}.
\end{align}
\textbf{Step 3.} Constraint $(\lambda-\lambda')(\mu-\mu')>0$:\\
The original constraint is equivalent to 
\[
    \frac{v}{u} = -\frac{ku-L}{ku+m_2AC}>0.
\]
Substitute $u = m_2C/\Delta$, and using 
\[
    ku-L = -m_1AC^2(|B|-B_1), \quad ku+m_2AC = AC(|B|-B_2), 
\]
where
\[
    B_1 = -\frac{m_2}{\lambda+2\mu}+\frac{m_2\tau}{\mu(\tau+1)}+\frac{m_2^2}{m_1\mu(\tau+1)}, \quad B_2=\frac{m_1}{\lambda+2\mu}+\frac{m_2(\tau-1)}{\mu(\tau+1)},
\]
We derive 
\begin{align}
    m_1m_2(|B|-B_1)(|B|-B_2)>0.
\end{align}
\textbf{Step 4.} Constraint $\lambda'+\mu'>0$:\\
We first consider the quantity $u+v$:
\[
u+v = u-\frac{u(ku-L)}{ku+m_2AC} = \frac{uSAC}{ku+m_2AC} = \frac{ACS}{k+A\Delta} = \frac{S}{|B|-B_2},
\]
where $S=m_1+m_2$. We notice that $\lambda'+\mu'=D+u+v$. Therefore, we have
\begin{align}
    D+\frac{S}{|B|-B_2}>0 \Rightarrow (|B|-B_2)\bra{|B|-B_2+\frac{m_1+m_2}{\lambda+\mu}}>0
\end{align}

\subsection{Proof of Proposition \ref{prop:volume}}\label{proof_elliptic_prop}
We use the notation introduced in the preceding proof. It remains only to derive the additional bounds resulting from the prescribed ranges of the Lam\'e parameters, since these assumptions imply the constraints used in Theorem \ref{thm:volume}.\\
\textbf{Step 1.} Constraint $\alpha^-\lambda<\lambda'<\alpha^+\lambda$:\\
By (C.6) and (C.7), we have
\begin{align}
    v = \frac{m_1m_2C(|B|-B_1)}{(|B|-B_2)|B|-B_3)},
\end{align}
where 
\[
    B_3 = \frac{m_2\tau}{C}.
\]
Then, we have the constraint,  
\[
    (\alpha^--1)\lambda<\frac{m_1m_2C(|B|-B_1)}{(|B|-B_2)(|B|-B_3)}<(\alpha^+-1)\lambda
\]
\textbf{Step 2.} Constraint $\beta^-\mu<\mu'<\beta^+\mu$:
Under the constraint $(\lambda-\lambda')(\mu-\mu')>0$, the admissible ranges for $\beta^\pm$ lead to the following two cases. \\
\textbf{Case 1.} $0<\beta^-<\beta^+<1$: We observe that $\mu'<\mu$, which implies $m_2<0$. Combining this observation and the fact that $u<0$, we derive that $\Delta>0$.\\
\textbf{Case 2.} $\beta^+>\beta^->1$: We observe that $\mu'>\mu$, which implies $m_2>0$. Combining this observation and the fact that $u>0$, we derive that $\Delta>0$.\\
Thus, we have
\[
    u = \frac{m_2C}{\Delta} \in \bra{(\beta^--1)\mu,(\beta^+-1)\mu},
\]
which implies 
\[
    (\beta^--1)\mu\bra{|B|C-m_2\tau}<m_2C<(\beta^+-1)\mu\bra{|B|C-m_2\tau}.
\]
Hence, we derive that
\begin{align}
    &|B|\in\bra{\frac{m_{2}}{\mu}\bra{\frac{\tau}{\tau+1}-\frac{1}{1-\beta^-}},\frac{m_{2}}{\mu}\bra{\frac{\tau}{\tau+1}-\frac{1}{1-\beta^+}}},\quad 0<\beta^-<\beta^+<1, \\
    &|B|\in\bra{\frac{m_{2}}{\mu}\bra{\frac{\tau}{\tau+1}-\frac{1}{1-\beta^+}},\frac{m_{2}}{\mu}\bra{\frac{\tau}{\tau+1}-\frac{1}{1-\beta^-}}},\quad \beta^+>\beta^->1
\end{align}

\section{Auxiliary integrals}
In this section, we present the following two integrals that are frequently used in the calculation of the point spread functions:
\begin{lem} For fixed $a>0$, $\varphi^s\in [0,2\pi)$, and $n\in \sZ$, we have
    \begin{align}
        &\int_0^{2\pi} \cos(n\theta)e^{ia\cos(\theta-\varphi^s)}d\theta= 2\pi i^n J_n(a)\cos(n\varphi^s),\\
        &\int_0^{2\pi} \sin(n\theta)e^{ia\cos(\theta-\varphi^s)}d\theta= 2\pi i^n J_n(a)\sin(n\varphi^s).
    \end{align}
\end{lem}
\begin{proof}
    We only provide the proof of the first integral; the second one can be proved similarly.
    By Jacobi–Anger expansion, we have 
    \[
        e^{ia\cos(\theta-\varphi^s)} = \sum_{m\in\sZ} i^mJ_m(a)e^{im(\theta-\varphi^s)}.
    \]
    Then, 
    \begin{align*}
        \int_0^{2\pi} \cos(n\theta)e^{ia\cos(\theta-\varphi^s)}d\theta &=  \sum_{m\in\sZ} i^mJ_m(a)\int_0^{2\pi} \cos(n\theta)e^{im(\theta-\varphi^s)} d\theta \\
        & = \pi i^n J_n(a)e^{-in\varphi^s}+\pi i^{-n} J_{-n}(a)e^{i n\varphi^s}\\
        & = 2\pi i^n J_n(a)\cos(n\varphi^s).
    \end{align*}
\end{proof}
Consequently, we obtain the following integrals.
\begin{lem}\label{lem:aux}
Let $\vr^s \in \sR^2\setminus\{\mathbf{0}\}$ be a vector, $r^s$ be its modulus and $a>0$ be a parameter. We introduce the polar coordinate $r^s_1 = r^s\cos\varphi^s$, $r^s_2 = r^s\sin\varphi^s$. For each $n\in\sZ$, we have
    \begin{align}
        &\int_0^{2\pi} \cos(n\theta)e^{ia(r^s_1\cos\theta+ r^s_2\sin\theta)}d\theta= 2\pi i^n J_n(ar^s)\cos(n\varphi^s),\\
        &\int_0^{2\pi} \sin(n\theta)e^{ia(r^s_1\cos\theta+ r^s_2\sin\theta)}d\theta= 2\pi i^n J_n(ar^s)\sin(n\varphi^s).
    \end{align}
\end{lem}
\begin{proof}
    We only need to notice that $r^s_1\cos\theta+ r^s_2\sin\theta = r^s\cos(\theta-\varphi^s)$.
\end{proof}

\newpage
\bibliography{reference}
\bibliographystyle{abbrv}
\end{document}